\documentclass[11pt]{article}
\usepackage{mathrsfs}

\usepackage{amssymb}
\usepackage{amsmath}
\usepackage{amsfonts}
\usepackage{graphics}
\usepackage{epsfig}
\usepackage{enumerate}

\usepackage{color}

\allowdisplaybreaks

\newtheorem{theorem}{Theorem}[section]
\newtheorem{lemma}[theorem]{Lemma}
\newtheorem{definition}[theorem]{Definition}
\newtheorem{proposition}[theorem]{Proposition}
\newtheorem{corollary}[theorem]{Corollary}
\newtheorem{remark}[theorem]{Remark}

\newenvironment{proof}{{\bf Proof:}}{~\hfill $\Box$}
\newenvironment{keywords}{{\bf Keywords: }}{}
\newenvironment{AMS}{{\bf AMS Subject Classification: }}{}

\numberwithin{equation}{section}

\begin{document}

\author{Qian Lin $^{1,2}$\footnote{ {\it Email address}: Qian.Lin@univ-brest.fr}
\\
{\small $^1$ School of Mathematics, Shandong University, Jinan 250100, People's Republic of China
} \\
{\small $^1$ Laboratoire de Math\'ematiques, CNRS UMR 6205,
Universit\'{e} de Bretagne Occidentale,}\\ {\small 6, avenue Victor
Le Gorgeu, CS 93837, 29238 Brest cedex 3, France.}  }

\title{Representation of G-martingales as stochastic integrals with respect to the G-Brownian motion
\thanks {This work is supported by the Marie Curie Initial
Training Network (PITN-GA-2008-213841).}}
\date{}

\maketitle \noindent
\begin{abstract}
The objective of this paper is to derive a representation of
symmetric G-martingales as stochastic integrals with respect to the G-Brownian
motion. For this end, we first study some extensions of stochastic
calculus with respect to G-martingales under the sublinear
expectation spaces.
\end{abstract}

\noindent
\begin{keywords}
G-expectation; G-Brownian motion; G-martingales; representation;
stochastic calculus.
\end{keywords}

\vskip 3mm \noindent
\begin{AMS}
    60H10; 60H05; 60H30.
\end{AMS}

\section{Introduction }

   Motivated by  uncertainty  problems, risk  measures  and
superhedging in finance,  Peng has introduced recently  a new notion
of nonlinear expectation, the so-called G-expectation (cf.
\cite{Peng:2006},  \cite{Peng:2007}, \cite{Peng:2008}),
   which, unlike the classical linear one, is not associated with the linear but a nonlinear heat
   equation.
The G-expectation represents a special case of general nonlinear
expectations $\mathbb{\hat{E}}$ which importance stems from the fact
that they are related to risk measures $\rho$ in finance by the
relation $\mathbb{\hat{E}}[X]=\rho(-X)$, where $X$ runs the class of
contingent claims. Although the G-expectations represent only a
special case, their importance inside the class of nonlinear
expectations reflects in the law of large numbers and central limit
theorem under nonlinear expectation, proven by Peng \cite{Peng2007}
and \cite{Peng2008-1}. Together with the notion of G-expectation
Peng also introduced the related G-normal distribution and the
G-Brownian motion. The $G$-Brownian motion is a stochastic process
which, under the $G$-expectation, has independent increments which
are G-normally distributed. Moreover, in \cite{Peng:2007}  Peng
developed an It\^{o} calculus for the G-Brownian motion.

A celebrated result of L\'{e}vy \cite{L1948} and  Doob \cite{Doob2}
states that a continuous classical martingale $M$ is a Brownian
motion if and only if its  quadratic variation process is the
deterministic function $\langle M\rangle_{t}=t, \ t\geq 0$.
Recently, the martingale characterization of the G-Brownian motion
has been obtained by Xu \cite{X:2009}.  The objective of the present
paper is to extend this characterization and to investigate a
representation of symmetric G-martingales as stochastic integrals with respect
to the G-Brownian motion in the framework of the sublinear
expectation spaces. In order to obtain this, we first study the
stochastic integral with respect to a larger class of symmetric G-martingales
$M$.  This generalizes considerably the recent works by Xu
\cite{X:2009}, in which the process $\{M^{2}_{t}-t\}_{t\in[0,T]}$
has been supposed to be a G-martingale.  We discuss even more general case in our another paper.

 However, there are  several difficulties for studying the martingale characterization of the G-Brownian
motion: firstly, in contrast to the classical Brownian motion the
G-Brownian motion is not defined on a given probability space but
only on a nonlinear expectation space. The nonlinear expectation
$\mathbb{\hat{E}}[\cdot]$ can be represented as the supremum over
the linear expectations $E_P[\cdot]$, where $P$ runs a certain class
of probability measures which are not mutually equivalent. Secondly,
the quadratic variation process $\langle B\rangle$ of the G-Brownian
motion is a random process which satisfies the relation
$\underline{\sigma}dt\le d\langle B\rangle_t\le
\overline{\sigma}dt,\ q.s., t\ge 0.$  Thirdly, related with their
absence  or restriction is the applicability of some well known
tools in the classical case (i.e. localization with stopping times,
the dominated convergence theorem).

Our paper is organized as follows: Section 2 introduces the
necessary notations and it gives a short recall of some elements of
the G-stochastic analysis, which will be used in what follows.
Moreover, the notion of G-martingales will be introduced. In Section
3  we extend  the stochastic calculus with respect to G-martingales
to a larger class of G-martingales and study related properties.
Moreover, a downcrossing inequality for G-supermartingales is
obtained. Finally,  section 4 investigates the representation of
G-martingales as stochastic integrals with respect to  G-Brownian
motion. This leads to the main result of this paper.

\section{Notations and preliminaries  }
In this section, we introduce the G-framework which was established
by Peng \cite{Peng:2006}, and which we will need in what follows.

Let $\Omega$ be a given nonempty set and $\mathcal {H}$ be a linear
space of real functions defined on $\Omega$ such that if
$x_{1},\cdot\cdot\cdot,x_{n}\in \mathcal {H}$ then
$\varphi(x_{1},\cdot\cdot\cdot,x_{n})\in \mathcal {H}$, for each
$\varphi\in C_{l,lip}(\mathbb{R}^{m})$. Here
$C_{l,lip}(\mathbb{R}^{m})$ denotes the linear space of functions
$\varphi$ satisfying
$$|\varphi(x)-\varphi(x)|\leq C(1+|x|^{n}+|y|^{n})|x-y|,\ \text{for all}\ x, y\in \mathbb{R}^{m},$$
for some $C>0$ and $\ n\in\mathbb{N}$, both depending on $\varphi$.
The space $\mathcal {H}$ is considered as a set of random variables.

\begin{definition}
{\bf A Sublinear expectation} $\mathbb{\hat{E}}$ on $\mathcal {H}$
is a functional $\mathbb{\hat{E}}:\mathcal {H}\mapsto \mathbb{R}$
satisfying the following properties: for all $X, Y \in \mathcal
{H}$, we have
\begin{enumerate}
\item[(i)] {\bf Monotonicity:} If $X\geq Y $, then
$\mathbb{\hat{E}}[X]\geq\mathbb{\hat{E}}[Y]$.
\item[(ii)] {\bf Constant preserving:}
$\mathbb{\hat{E}}[c]=c$, for all $c\in \mathbb{R}$.
\item[(iii)] {\bf Self-dominated property:}
$\mathbb{\hat{E}}[X]-\mathbb{\hat{E}}[Y]\leq\mathbb{\hat{E}}[X-Y]$.
\item[(iv)] {\bf Positive homogeneity:}
$\mathbb{\hat{E}}[\lambda X]=\lambda\mathbb{\hat{E}}[X],$ for all $
\lambda \geq 0$.
\end{enumerate}
The triple $(\Omega, \mathcal {H}, \mathbb{\hat{E}})$ is called a
sublinear expectation space.
\end{definition}

\begin{remark}
    The sublinear expectation space can be regarded as a generalization of the classical
probability space $(\Omega, \mathcal {F}, \mathbb{P})$ endowed with
the linear expectation associated with $\mathbb{P}$.
\end{remark}

\begin{definition} In a sublinear expectation space $(\Omega, \mathcal {H},
\mathbb{\hat{E}})$, a random vector $Y=(Y_{1},\cdots, Y_{n}),
Y_{i}\in \mathcal {H}$, is said to be independent under $
\mathbb{\hat{E}}$ of another random vector $X=(X_{1},\cdots, X_{m}),
X_{i}\in \mathcal {H}$,  if for each test function $\varphi \in
C_{l,lip}(\mathbb{R}^{m+n})$ we have
$$\mathbb{\hat{E}}[\varphi(X, Y)]=\mathbb{\hat{E}}[\mathbb{\hat{E}}[\varphi(x, Y)]_{x=X}].$$
\end{definition}

\begin{definition}{\bf(G-normal distribution)}
 Let be given two reals $\underline{\sigma}, \overline{\sigma}$ with
$0\leq\underline{\sigma}\leq \overline{\sigma}.$ A random variable
$\xi$ in a sublinear expectation space $(\Omega, \mathcal {H},
\mathbb{\hat{E}})$ is called $G_{\underline{\sigma},
\overline{\sigma}}$-normal distributed, denoted by $\xi\sim \mathcal
{N}(0, [{\underline{\sigma}^{2}, \overline{\sigma}}^{2}])$,  if for
each $\varphi \in C_{l,lip}(\mathbb{R})$, the following function
defined by
$$u(t, x):=\mathbb{\hat{E}}[\varphi(x+\sqrt{t}\xi)], \quad (t,x)\in [0, \infty)\times\mathbb{R},$$
is the unique continuous  viscosity solution with polynomial growth
of the following parabolic partial differential equation :
$$\left\{\begin{array}{l l}
      \partial_{t}u(t, x)=G(\partial_{xx}^{2}u(t, x)),\quad (t, x)\in[0,
\infty)\times\mathbb{R},\\
       u(0, x)=\varphi(x).
         \end{array}
  \right.$$
 Here  $G=G_{\underline{\sigma},
\overline{\sigma}}$
 is the following sublinear function parameterized by $\underline{\sigma}$
 and $\overline{\sigma}$:
 $$G(\alpha)=\frac{1}{2}(\overline{\sigma}^{2}\alpha^{+}-\underline{\sigma}^{2}\alpha^{-}), \quad \alpha \in \mathbb{R}$$
 (Recall that $\alpha^{+}=\text{max}\{0,\alpha\}$ and
$\alpha^{-}=-\text{min}\{0,\alpha\}$).
\end{definition}
For simplicity, we suppose that  $\overline{\sigma}^{2}=1$ and $
\underline{\sigma}^{2}=\sigma^{2}_{0}, \sigma^{2}_{0}\in [0, 1]$, in
the following paper.

 Throughout this paper, we let $\Omega=C_{0}(\mathbb{R^{+}})$ be the space of all
real valued continuous functions $(\omega_{t})_{t\in
\mathbb{R^{+}}}$ with $\omega_{0}=0$, equipped with the distance
$$\rho(\omega^{1}, \omega^{2})
=\sum\limits_{i=1}^{\infty}2^{-i}\Big[(\max\limits_{t\in[0,i]}|\omega_{t}^{1}-\omega_{t}^{2}|)\wedge1\Big],
\ \omega^{1}, \omega^{2}\in\Omega.$$

For each $T > 0$, we consider the following space of random
variables:
\begin{eqnarray*}
L_{ip}^{0}(\mathcal
{F}_{T}):=\Big\{X(\omega)=\varphi(\omega_{{t_{1}}} \cdots,
\omega_{{t_{m}}}) \ | \ t_{1}, \cdots, t_{m}\in [0,T],  \text{ for
all} \ \varphi\in C_{l,lip}(\mathbb{R}^{m}), \ m\geq1\Big\}.
\end{eqnarray*}
Obviously, it holds that $L_{ip}^{0}(\mathcal {F}_{t}) \subseteq
L_{ip}^{0}(\mathcal {F}_{T})$, for all $t\leq T<\infty$. We notice
that $X,Y \in L_{ip}^{0}(\mathcal {F}_{t})$ implies $X\cdot Y \in
L_{ip}^{0}(\mathcal {F}_{t})$ and $|X|\in L_{ip}^{0}(\mathcal
{F}_{t})$. We further define
$$L_{ip}^{0}(\mathcal {F})=\bigcup_{n=1}^{\infty}L_{ip}^{0}(\mathcal {F}_{n}).$$

We will work on the canonical space $\Omega$ and set
$B_{t}(\omega)=\omega_{t}$, $t\in [0,\infty)$, for $\omega\in
\Omega.$

We now introduce a sublinear expectation $\mathbb{\hat{E}}$ defined
 on $\mathcal {H}_{T}^{0}=L_{ip}^{0}(\mathcal
{F}_{T})$ as well as on $\mathcal {H}^{0}=L_{ip}^{0}(\mathcal {F})$.
For this, we consider the function $G(a)= \frac{1 }{2}
(a^{+}-\sigma^{2}_{0}a^{-}), a \in \mathbb{R}$, and we apply the
following procedure: for each $X\in \mathcal {H}^{0}$ with
$$X=\varphi(B_{t_{1}}-B_{t_{0}}, B_{t_{2}}-B_{t_{1}}, \cdots, B_{t_{m}}-B_{t_{m-1}})$$
for some $m\geq 1, \varphi\in C_{l,lip}(\mathbb{R}^{m})$ and
$0=t_{0}\leq t_{1}\leq\cdots \leq t_{m}<\infty$, we set
\begin{eqnarray*}
&&\mathbb{\hat{E}}[\varphi(B_{t_{1}}-B_{t_{0}},
B_{t_{2}}-B_{t_{1}}, \cdots, B_{t_{m}}-B_{t_{m-1}})]\\
&&=\mathbb{\widetilde{E}}[\varphi(\sqrt{t_{1}-t_{0}}\xi_{1},
\sqrt{t_{2}-t_{1}}\xi_{2}, \cdots, \sqrt{t_{m}-t_{m-1}}\xi_{m})],
\end{eqnarray*}
where $(\xi_{1},\xi_{2},\cdots, \xi_{m})$ is an m-dimensional
G-normal distributed random vector in a sublinear expectation space
$(\widetilde{\Omega}, \widetilde{\mathcal {H}},
\widetilde{\mathbb{E}})$ such that $\xi_{i}\sim\mathcal {N}(0,
[\sigma^{2}_{0},1])$ and  $\xi_{i+1}$ is independent of $(\xi_{1},
\cdots, \xi_{i}),$ for every $i=1, 2, \cdots, m$.

The related conditional expectation of
$X=\varphi(B_{t_{1}}-B_{t_{0}}, B_{t_{2}}-B_{t_{1}}, \cdots,
B_{t_{m}}-B_{t_{m-1}})$ under $\mathcal {H}_{t_{j}}^{0}$ is defined
by
\begin{eqnarray*}
\mathbb{\hat{E}}[X|\mathcal
{H}_{t_{j}}^{0}]&=&\mathbb{\hat{E}}[\varphi(B_{t_{1}}-B_{t_{0}},
B_{t_{2}}-B_{t_{1}}, \cdots, B_{t_{m}}-B_{t_{m-1}})|\mathcal
{H}_{t_{j}}^{0}]\\
&=&\psi(B_{t_{1}}-B_{t_{0}}, B_{t_{2}}-B_{t_{1}}, \cdots,
B_{t_{j}}-B_{t_{j-1}}),
\end{eqnarray*}
where
$$\psi(x_{1}, x_{2}, \cdots,
x_{j})=\mathbb{\widetilde{E}}[\varphi(x_{1}, x_{2}, \cdots, x_{j},
\sqrt{t_{j+1}-t_{j}}\xi_{j+1}, \cdots,
\sqrt{t_{m}-t_{m-1}}\xi_{m})],$$ $(x_{1}, x_{2}, \cdots,
x_{j})\in\mathbb{R}^{j}, 0\leq j\leq m.$

\vspace{2mm} For $p\geq 1$,
$\|X\|_{p}=\mathbb{\hat{E}}^{\frac{1}{p}}[|X|^{p}]$,  $X \in
L_{ip}^{0}(\mathcal {F}),$ defines a norm on  $L_{ip}^{0}(\mathcal
{F})$. Let $\mathcal {H}=L_{G}^{p}(\mathcal {F})$ (resp. $\mathcal
{H}_{t}=L_{G}^{p}(\mathcal {F}_{t})$) be the completion of
$L_{ip}^{0}(\mathcal {F})$ (resp. $L_{ip}^{0}(\mathcal {F}_{t})$)
under the norm $\|\cdot\|_{p}$. Then the space $(L_{G}^{p}(\mathcal
{F}), \|\cdot\|_{p})$ is a Banach space and the operators
$\mathbb{\hat{E}}[\cdot] $ and $\mathbb{\hat{E}}[\cdot|\mathcal
{H}_{t}]$ can be continuously extended to the Banach space
$L_{G}^{p}(\mathcal {F})$. Moreover, we have $L_{G}^{p}(\mathcal
{F}_{t})\subseteq L_{G}^{p}(\mathcal {F}_{T})\subset
L_{G}^{p}(\mathcal {F})$, for all $0\leq t \leq T <\infty$.

\begin{definition}
  The expectation $\mathbb{\hat{E}}: L_{G}^{p}(\mathcal {F})\mapsto \mathbb{R}$  defined through the above procedure is
called G-expectation.
\end{definition}

\begin{proposition}\label{p2}
  For all
$t, s\in [0, \infty)$, we list the properties of
$\mathbb{\hat{E}}[\cdot|\mathcal {H}_{t}]$ that hold for all $X, Y
\in L_{G}^{p}(\mathcal {F}):$

\begin{enumerate}[(i)]
\item  If $X\geq Y$, then $\mathbb{\hat{E}}[X|\mathcal
{H}_{t}]\geq \mathbb{\hat{E}}[Y|\mathcal {H}_{t}]$;
\item $\mathbb{\hat{E}}[\eta|\mathcal
{H}_{t}]=\eta$, for all $\eta\in L_{G}^{p}(\mathcal {F}_{t})$;
\item $\mathbb{\hat{E}}[X|\mathcal
{H}_{t}]- \mathbb{\hat{E}}[Y|\mathcal
{H}_{t}]\leq\mathbb{\hat{E}}[X-Y|\mathcal {H}_{t}];$
\item $\mathbb{\hat{E}}[\eta X|\mathcal
{H}_{t}]=\eta^{+}\mathbb{\hat{E}}[X|\mathcal
{H}_{t}]-\eta^{-}\mathbb{\hat{E}}[- X|\mathcal {H}_{t}]$, for all
$\eta\in L_{G}^{p}(\mathcal {F}_{t})$;
\item If $\mathbb{\hat{E}}[Y|\mathcal
{H}_{t}]=-\mathbb{\hat{E}}[-Y|\mathcal {H}_{t}],$  then
$\mathbb{\hat{E}}[X+Y|\mathcal {H}_{t}]= \mathbb{\hat{E}}[X|\mathcal
{H}_{t}]+\mathbb{\hat{E}}[Y|\mathcal {H}_{t}];$
\item $\mathbb{\hat{E}}[\mathbb{\hat{E}}[X|\mathcal
{H}_{t}]|\mathcal {H}_{s}]=\mathbb{\hat{E}}[X|\mathcal {H}_{t\wedge
s}]$, and, in particular,
$\mathbb{\hat{E}}[\mathbb{\hat{E}}[X|\mathcal
{H}_{t}]]=\mathbb{\hat{E}}[X]$.
\end{enumerate}

\end{proposition}

 For $p\geq 1$ and an arbitrary but fixed  time horizon $0<T<\infty$, we now consider the
following space of step processes:
\begin{eqnarray*}
M_{G}^{p, 0}(0,T) &=& \bigg\{\eta:
\eta_{t}=\sum\limits_{j=0}^{n-1}\xi_{j}I_{[t_{j},t_{j+1})},
0=t_{0}<t_{1}<\cdots<t_{n}=T, \\&& \hskip 1cm\xi_{j}\in
L_{G}^{p}(\mathcal {F}_{t_{j}}), j=0,\cdots, n-1, \text{for all}\
n\geq 1\bigg\},
\end{eqnarray*}
and we  define the following norm in $M_{G}^{p, 0}(0,T)$:

$$\parallel \eta
\parallel_{p}=\bigg(\mathbb{\hat{E}}\Big[\int_{0}^{T}|\eta_{t}|^{p}dt\Big]\bigg)^{\frac{1}{p}}=\bigg(\mathbb{\hat{E}}
\Big[\sum\limits_{j=0}^{n-1}|\xi_{t_{j}}|^{p}(t_{j+1}-t_{j})\Big]\bigg)^{\frac{1}{p}}.
$$
Finally, we denote by $M_{G}^{p}(0,T)$ the completion of $M_{G}^{p,
0}(0,T)$ under the norm $\parallel \cdot\parallel_{p}.$

\begin{definition}
 A process $B=\{B_{t},t\geq 0\}$ in a sublinear expectation space $(\Omega, \mathcal {H},
 \mathbb{\hat{E}})$ is called a G-Brownian motion if $ \{B_{t}, t\geq 0\}\subset\mathcal {H}$ and the following
 properties are satisfied:
\begin{enumerate}
\item[(i)] $B_{0}=0$;
\item[(ii)] For each $t, s \geq0$, the difference $B_{t+s}-B_{t}$ is
$\mathcal {N}(0, \ [\sigma_{0}^{2}s, s])$-distributed and is
independent of $(B_{t_{1}}, \cdots, B_{t_{n}})$, for all
$n\in\mathbb{N}$ and $0\leq t_{1}\leq\cdots \leq t_{n}\leq t$.
\end{enumerate}
\end{definition}

\begin{remark}
The  canonical process $(B_{t})_{t\geq 0}$ in $(\Omega, \mathcal
{H})$, $\Omega=C_{0}(\mathbb{R}_{+})$, endowed with the
G-expectation $\mathbb{\hat{E}}$ is a G-Brownian motion.

\end{remark}

\begin{remark}
  In \cite{Peng:2006}, \cite{Peng:2007} and \cite{Peng:2008}, Peng established a stochastic calculus of
It\^o's type with respect to the G-Brownian motion  and its
quadratic variation process.  Peng derived an It\^o's formula and
moreover, he obtained the existence and uniqueness of the solution
to stochastic differential equations with Lipschitz coeffcients
driven by G-Brownian motion.
\end{remark}

In \cite{HP:2009}, Hu and Peng obtained the representation theorem
of G-Expectations as follows.
\begin{proposition}
 Let $\mathbb{\hat{E}}$ be G-expectation. Then there exists a weekly compact
family of probability measures $\mathcal {P}$ on $(\Omega, \mathcal
{B}(\Omega))$ such that
$$\mathbb{\hat{E}}[X]=\max\limits_{P\in\mathcal {P}}E_{P}[X], \ \text{for
all}\ X\in \mathcal {H},$$ where $E_{P}[\cdot]$ is the linear
expectation with respect to $P$.
\end{proposition}
The authors of \cite{HP:2009} also introduced the associate Choquet
capacity $$c(A)=\sup\limits_{P\in\mathcal {P}}P(A), A\in\mathcal
{B}(\Omega). $$ We have the following proposition.
\begin{proposition}\label{pro}
\begin{enumerate}
\item[(i)]   $0\leq c(A)\leq 1 $, for all $A\subset\Omega$.
\item[(ii)] If $A\subset B$, then $c(A)\leq c(B)$.
\item[(iii)] If $\{A_{n}\}_{n=1}^{\infty}$ is an increasing sequence
in $\mathcal{B}(\Omega)$ and  $A_{n}\uparrow A$, then
$c(A)=\lim\limits_{n\rightarrow\infty}c(A_{n})$.
\end{enumerate}
\end{proposition}

 \begin{definition}
 A set $A$ is polar if $c(A)=0$ and a property holds quasi-surely (q.s.) if it holds
 outside a polar set.
\end{definition}

  As in the classical stochastic analysis, the definition of a
  modification of a process plays an important role.
\begin{definition}
  Let $I$ be a set of indexes, and $\{X_{t}\}_{t\in I}$ and $\{Y_{t}\}_{t\in
  I}$ two processes indexed by $I$. We say that $Y$ is a
  modification of $X$ if for all $t\in I$, $X_{t}=Y_{t}\ q.s.$
\end{definition}

Finally, we recall the definition of a G-martingale introduced by
Peng \cite{Peng:2008}.

\begin{definition}
 A process $M=\{M_{t},t\geq 0\}$ is called a G-martingale (respectively,
 G-supermartingale, and G-submartingale) if for each  $t\in [0, \infty),
 M_{t}\in L_{G}^{1}(\mathcal {F}_{t}) $  and for each  $s\in [0, t]$, we
 have
 $$\mathbb{\hat{E}}[M_{t}|\mathcal {H}_{s}]=M_{s}, \ (\text{respectively}\ \leq M_{s}, \text{and}\  \geq
 M_{s})\ q.s. $$
A process $M=\{M_{t},t\geq 0\}$ is called a symmetric G-martingale,
if $M$ and $-M$ are G-martingales.
\end{definition}

\section{Stochastic integrals of G-martingales}

In this sections, we study the stochastic integrals of G-martingales
and related properties, which will be important in next section.

Let $p\geq 1$ and  $T>0$ be an arbitrarily fixed time horizon. Let
$\{A_{t},t\in [0,T]\}$ be a continuous and increasing process such
that for all $t\in[0,T], A_{t}\in \mathcal {H}_{t}$, $ A_{0}=0$ and
$\mathbb{\hat{E}}[A_{T}]<\infty$.
 We first consider the
following space of step processes:
\begin{eqnarray*}
M_{G}^{p, 0}(0,T) &=& \bigg\{\eta:
\eta_{t}=\sum\limits_{j=0}^{n-1}\xi_{t_{j}}I_{[t_{j},t_{j+1})},
0=t_{0}<t_{1}<\cdots<t_{n}=T, \\&& \hskip 1cm\xi_{t_{j}}\in
L_{G}^{p}(\mathcal {F}_{t_{j}}), j=0,\cdots, n-1, \text{for all}\
n\geq 1\bigg\},
\end{eqnarray*}
and we  define the following norm in $M_{G}^{p, 0}(0,T)$:

$$\parallel \eta \parallel_{p}=\bigg(\mathbb{\hat{E}}
\Big[\int_{0}^{T}|\eta_{t}|^{p}dA_{t}\Big]\bigg)^{\frac{1}{p}}
=\bigg(\mathbb{\hat{E}}
\Big[\sum\limits_{j=0}^{n-1}|\xi_{t_{j}}|^{p}(A_{t_{j+1}}-A_{t_{j}})\Big]\bigg)^{\frac{1}{p}}.$$
We denote by $M_{G,A}^{p}(0,T)$ the completion of $M_{G}^{p,
0}(0,T)$ under the norm $\parallel \cdot\parallel_{p}$, and  we
introduce the following space of G-martingales related with $A$:
$$\mathcal {M}=\Big\{M| M \ \text{is a continuous  symmetric
G-martingale such that} \ M^{2}-A\  \ \text{is a
G-supermartingale}\Big\}.$$

We will see later  that $\mathcal {M}\subset M_{G,A}^{2}(0,T)$.
\begin{definition}
For any $M\in\mathcal {M}$ and $\eta\in M_{G}^{2,0}(0,T)$ of the
form
 $\eta_{t}=\sum\limits_{j=0}^{n-1}\xi_{t_{j}}I_{[t_{j},t_{j+1})}(t),$
 we define $$I(\eta)=\int_0^T
 \eta_{t}dM_{t}=\sum\limits_{j=0}^{n-1}\xi_{t_{j}}(M_{t_{j+1}}
 -M_{t_{j}}).$$
\end{definition}

\begin{proposition}\label{pr1}
For all $M\in\mathcal {M}$, the mapping $I:
M_{G}^{2,0}(0,T)\rightarrow L_{G}^{2}(\mathcal {F}_{T})$ is a linear
continuous mapping and thus can be continuously extended to $I:
M_{G,A}^{2}(0,T)\rightarrow L_{G}^{2}(\mathcal {F}_{T})$. Moreover,
for all $\eta\in M_{G,A}^{2}(0,T)$, the process $\Big\{\int_0^t
 \eta_{s}dM_{s}\Big\}_{t\in[0,T]}$ is a symmetric G-martingale and
\begin{eqnarray}\label{ineq1}
\mathbb{\hat{E}}\Big[|\int_0^T \eta_{t}dM_{t}|^{2}\Big]
 \leq \mathbb{\hat{E}}\Big[\int_0^T |\eta_{t}|^{2}dA_{t}\Big].
\end{eqnarray}
\end{proposition}

\begin{proof}
From $M$ is a symmetric G-martingale and $M^{2}-A\ $ is a
G-supermartingale it follows that, for all $0\leq s \leq t \leq T$,
\begin{eqnarray*}
&&\mathbb{\hat{E}}[(M_{t}-M_{s})^{2}-(A_{t}-A_{s})|\mathcal {H}_{s}]\\
&=&\mathbb{\hat{E}}[M_{t}^{2}-M_{s}^{2}-2M_{s}(M_{t}-M_{s})-(A_{t}-A_{s})|\mathcal
{H}_{s}]\\
&=&\mathbb{\hat{E}}[M_{t}^{2}-M_{s}^{2}-(A_{t}-A_{s})|\mathcal
{H}_{s}]\leq 0.
\end{eqnarray*}

For $\eta\in M_{G}^{2,0}(0,T)$ of the form
 $\eta_{t}=\sum\limits_{j=0}^{n-1}\xi_{t_{j}}I_{[t_{j},t_{j+1})}(t),$ we
 have
\begin{eqnarray*}
 \mathbb{\hat{E}}\Big[|\int_0^T \eta_{t}dM_{t}|^{2}\Big]
 &= &\mathbb{\hat{E}}\Big[|\sum\limits_{j=0}^{n-1}\xi_{t_{j}}(M_{t_{j+1}}
 -M_{t_{j}})|^{2}\Big]
=\mathbb{\hat{E}}\Big[\sum\limits_{j=0}^{n-1}\xi_{t_{j}}^{2}(M_{t_{j+1}}
 -M_{t_{j}})^{2}\Big]\\
 &=&\mathbb{\hat{E}}\Big[\sum\limits_{j=0}^{n-1}\xi_{t_{j}}^{2}(M_{t_{j+1}}^{2}
 -M_{t_{j}}^{2})\Big]\\
  &\leq& \mathbb{\hat{E}}\Big[\sum\limits_{j=0}^{n-1}\xi_{t_{j}}^{2}(M_{t_{j+1}}^{2}
 -M_{t_{j}}^{2}-A_{t_{j+1}}
 +A_{t_{j}})\Big]+\mathbb{\hat{E}}\Big[\sum\limits_{j=0}^{n-1}\xi_{t_{j}}^{2}(A_{t_{j+1}}
 -A_{t_{j}})\Big],
 \end{eqnarray*}
 where
 \begin{eqnarray*}
 &&\mathbb{\hat{E}}\Big[\sum\limits_{j=0}^{n-1}\xi_{t_{j}}^{2}(M_{t_{j+1}}^{2}
 -M_{t_{j}}^{2}-A_{t_{j+1}}
 +A_{t_{j}})\Big]\\ &\leq& \mathbb{\hat{E}}\Big[\sum\limits_{j=0}^{n-2}\xi_{t_{j}}^{2}(M_{t_{j+1}}^{2}
 -M_{t_{j}}^{2}-A_{t_{j+1}}
+A_{t_{j}}) + \\&&\xi_{t_{n-1} }^{2} \mathbb{\hat{E}}[(M_{t_{n}}^{2}
 -M_{t_{n-1}}^{2}-A_{t_{n}}
 +A_{t_{n-1}}|\mathcal {H}_{t_{n-1}}]\Big]\\
 &\leq& \mathbb{\hat{E}}\Big[\sum\limits_{j=0}^{n-2}\xi_{t_{j}}^{2}(M_{t_{j+1}}^{2}
 -M_{t_{j}}^{2}-A_{t_{j+1}}
+A_{t_{j}})\Big]\leq\cdots\leq 0.
 \end{eqnarray*}
 Consequently,
 \begin{eqnarray*}
 \mathbb{\hat{E}}\Big[|\int_0^T \eta_{t}dM_{t}|^{2}\Big]
 &\leq& \mathbb{\hat{E}}\Big[\sum\limits_{j=0}^{n-1}\xi_{t_{j}}^{2}(A_{t_{j+1}}
 -A_{t_{j}})\Big]=\mathbb{\hat{E}}\Big[\int_0^T
 |\eta_{t}|^{2}dA_{t}\Big].
 \end{eqnarray*}
 Thus, (\ref{ineq1}) holds for all $\eta\in M_{G}^{2,0}(0,T)$. We
then can continuously extend the above inequality to the case
$\eta\in M_{G,A}^{2}(0,T)$ and obtain (\ref{ineq1}).

For $\eta\in M_{G,A}^{2}(0,T)$, there exists  a sequence of
$\eta^{n}\in M_{G}^{2,0}(0,T)$ of the form
 $\eta_{t}^{n}=\sum\limits_{j=0}^{n-1}\xi_{t_{j}}I_{[t_{j},t_{j+1})}(t),$ $ \xi_{t_{j}}\in
L_{G}^{2}(\mathcal {F}_{t_{j}})$
 such that $$\mathbb{\hat{E}}[|\int_0^t
(\eta_{u}-\eta_{u}^{n})dM_{u}|^{2}]\rightarrow 0,\ \text{as}\
n\rightarrow\infty.$$

 Let $0\leq s \leq t \leq T$. Without loss of generality, we
 assume that $t_{i}\leq s < t_{i+1} < t,$ for some $0\leq i \leq n-1$. Then we have
\begin{eqnarray*}
\mathbb{\hat{E}}[\int_0^t \eta_{u}^{N}dM_{u}|\mathcal
{H}_{s}]&=&\mathbb{\hat{E}}[\sum\limits_{j=0}^{n-1}\xi_{t_{j}}(M_{t_{j+1}\wedge
t} -M_{t_{j}\wedge t)}|\mathcal{H}_{s}]\\
&=&\sum\limits_{j=0}^{i-1}\xi_{t_{j}}(M_{t_{j+1}}
-M_{t_{j}})+\xi_{t_{i}}(M_{s} -M_{t_{i}})=\int_0^s
\eta_{u}^{n}dM_{u}.
\end{eqnarray*}
Consequently, $\int_0^\cdot \eta_{u}^{n}dM_{u}$ is a G-martingale.
Moreover,
\begin{eqnarray*}
&&\mathbb{\hat{E}}[|\mathbb{\hat{E}}[\int_0^t
\eta_{u}dM_{u}|\mathcal {H}_{s}]-\int_0^s
\eta_{u}dM_{u}|^{2}]\\&=&\mathbb{\hat{E}}[|\mathbb{\hat{E}}[\int_0^t
\eta_{u}dM_{u}|\mathcal {H}_{s}]-\int_0^s
\eta_{u}dM_{u}-\mathbb{\hat{E}}[\int_0^t \eta_{u}^{n}dM_{u}|\mathcal
{H}_{s}]+\int_0^s\eta_{u}^{n}dM_{u}|^{2}]\\
&\leq&2\mathbb{\hat{E}}[|\int_0^t
(\eta_{u}-\eta_{u}^{n})dM_{u}|^{2}]+2\mathbb{\hat{E}}[|\int_0^t
(\eta_{u}-\eta_{u}^{n})dM_{u}|^{2}]\\
&\rightarrow &0,\ \text{as}\ N\rightarrow\infty.
\end{eqnarray*}
Therefore $$\mathbb{\hat{E}}[\int_0^t \eta_{u}dM_{u}|\mathcal
{H}_{s}]=\int_0^s \eta_{u}dM_{u},\ q.s.$$ which means that
$\int_0^\cdot \eta_{u}dM_{u}$ is a G-martingale. On the other hand,
by extending the associate property of $\int_0^\cdot
\eta_{u}^{n}dM_{u}$, we have
$\int_0^\cdot(-\eta_{u})dM_{u}=-\int_0^\cdot \eta_{u}dM_{u}$, so
that
$$\mathbb{\hat{E}}[-\int_0^t \eta_{u}dM_{u}|\mathcal
{H}_{s}]=-\int_0^s \eta_{u}dM_{u},\ q.s., 0\leq s \leq t.$$ Thus,
$$\mathbb{\hat{E}}[-\int_0^t \eta_{u}dM_{u}|\mathcal
{H}_{s}]=-\mathbb{\hat{E}}[\int_0^t \eta_{u}dM_{u}|\mathcal
{H}_{s}]=-\int_0^s \eta_{u}dM_{u},\ q.s.$$ Consequently, $\{\int_0^t
 \eta_{s}dM_{s}, t\in[0,T]\}$  is a symmetric G-martingale. The proof is complete.
\end{proof}

 For $0\leq s \leq t \leq T$ and $\eta\in M_{G,A}^{2}(0,T)$, we denote

\begin{eqnarray*}
\int_s^t \eta_{u}dM_{u}=\int_0^T I_{[s,t]}(u)\eta_{u}dM_{u}.
\end{eqnarray*}
It is now straightforward to see that we have the following
properties of the stochastic integral of G-martingales.

\begin{proposition}\label{pr5}
Let $0\leq s < r \leq t \leq T$. For all  $M\in\mathcal {M}$ and
$\theta, \eta\in M_{G,A}^{2}(0,T)$, we have

\begin{enumerate}[(i)]
\item   $\int_s^t \eta_{u}dM_{u}=\int_s^r \eta_{u}dM_{u}+\int_r^t \eta_{u}dM_{u}$;
\item $\int_s^t (\eta_{u}+\alpha\theta_{u}) dM_{u}=\int_s^t \eta_{u}dM_{u}+\alpha\int_s^t \theta_{u}dM_{u}$,
for all $\alpha$  bounded random variable in $L_{G}^{p}(\mathcal
{F}_{s})$;
\item $\mathbb{\hat{E}}[X+\int_r^T
\eta_{u}dM_{u}|\mathcal {H}_{s}]$=$\mathbb{\hat{E}}[X|\mathcal
{H}_{s}]$, for all $X\in L_{G}^{p}(\mathcal {F}).$
\end{enumerate}
\end{proposition}

For proving the continuity of the stochastic integral regarded as a
process, we need the following Doob inequality for symmetric
G-martingale.

\begin{theorem}\label{th1}
 If $X$ is a right-continuous symmetric G-martingale running over an interval $[0,T]$
 of  $\mathbb{R}$,  then for every $p>1$ such that $X_{T}\in L_{G}^{p}(\mathcal {F})$,
\begin{eqnarray*}
\mathbb{\hat{E}}[\sup\limits_{0\leq t\leq T} |X_{t}|^{p}]\leq
(\frac{p}{p-1})^{p}\mathbb{\hat{E}}[ |X_{T}|^{p}].
\end{eqnarray*}
\end{theorem}

\begin{proof}
By Remark \ref{re1}, there exists a weekly compact family of
probability measures $\mathcal {P}$ on $(\Omega, \mathcal
{B}(\Omega))$ such that
$\mathbb{\hat{E}}[X]=\max\limits_{P\in\mathcal {P}}E_{P}[X]$, for
all $X\in \mathcal {H}$, where $E_{P}[\cdot]$ is the linear
expectation with respect to $P$.

For all $0 \leq t\leq T$, let $\mathcal {F}_{t}^{B}=\sigma\{B_{s},
s\leq t\}$.  For any $0\leq s \leq t\leq T$ and any positive $\xi\in
L_{G}^{\frac{p}{p-1}}(\mathcal {F}_{s})$, we have
\begin{eqnarray*}
\mathbb{\hat{E}}[ (X_{t}-X_{s})\xi]=\mathbb{\hat{E}}[
\xi(\mathbb{\hat{E}}[X_{t}|\mathcal {H}_{s}]-X_{s})]= 0.
\end{eqnarray*}
On the other hand,
\begin{eqnarray*}
\mathbb{\hat{E}}[ (X_{t}-X_{s})\xi]=\max\limits_{P\in\mathcal
{P}}E_{P}[(X_{t}-X_{s})\xi] \geq E_{P}[(X_{t}-X_{s})\xi]= E_{P}[
\xi(E_{P}[X_{t}|\mathcal {F}_{s}^{B}]-X_{s})],
\end{eqnarray*}
then we have $E_{P}[X_{t}|\mathcal {F}_{s}^{B}]\leq X_{s}, P$-$
a.s,$ for all $P\in\mathcal {P}$. By the same argument but this time
with negative $\xi\in L_{G}^{\frac{p}{p-1}}(\mathcal {F}_{s})$, we
can prove that $E_{P}[X_{t}|\mathcal {F}_{s}^{B}]\geq X_{s}, P$-$
a.s,$ for all $P\in\mathcal {P}$. Therefore $E_{P}[X_{t}|\mathcal
{F}_{s}^{B}]= X_{s}, P$-$ a.s.,$  for all $P\in\mathcal {P}$. Thus
$X$ is a $P$-martingale and from the classical Doob's inequality it
follows that
\begin{eqnarray*}
E_{P}[\sup\limits_{0\leq t\leq T} |X_{t}|^{p}]\leq
(\frac{p}{p-1})^{p}E_{P}[ |X_{T}|^{p}], \ \text{for all}\
P\in\mathcal {P}.
\end{eqnarray*}
Therefore,
\begin{eqnarray*}
\mathbb{\hat{E}}[\sup\limits_{0\leq t\leq T} |X_{t}|^{p}]\leq
(\frac{p}{p-1})^{p}\mathbb{\hat{E}}[ |X_{T}|^{p}].
\end{eqnarray*} The proof is complete.
\end{proof}

We now give a downcrossing  inequality for G-supermartingales. Let
$a, b$ be two positive constants such that $a<b$. Let
$\pi_{n}=\{0=t_{0}<\cdots < t_{n}=T\}$ be a partition of the
interval $[0,T]$. We define $D_{a}^{b}[X,n]$  the number of
downcrossing of $[a,b]$ by $\{X_{t_{i}}\}_{i=0}^{n}$.
\begin{theorem}\label{}
Let $X$ be a positive G-supermartingale and $0=t_{0}\leq\cdots \leq
t_{n}=T$ be a strictly increasing sequences. Then for all real
positive numbers $a$ and $ b$ such that $a<b$,
\begin{eqnarray*}
\mathbb{\hat{E}}[D_{a}^{b}[X,n]]\leq \frac{1}{b-a}\mathbb{\hat{E}}[
X_{0}\wedge b].
\end{eqnarray*}
\end{theorem}

\begin{proof}
 For any $0\leq s
\leq t\leq T$, from the first part of  the proof of Theorem
\ref{th1}, we know that $E_{P}[X_{t}|\mathcal {F}_{s}^{B}]\leq
X_{s}, P$-$ a.s.,$ for all $P\in\mathcal {P}$.
 From the classical downcrossing  inequality for supermartingales (cf. \cite{Doob}) it follows that
\begin{eqnarray*}
E_{P}[D_{a}^{b}[X,n]]\leq \frac{1}{b-a}E_{P}[ X_{0}\wedge b], \
\text{for all}\ P\in\mathcal {P}.
\end{eqnarray*}
Therefore,
\begin{eqnarray*}
\mathbb{\hat{E}}[D_{a}^{b}[X,n]]\leq \frac{1}{b-a}\mathbb{\hat{E}}[
X_{0}\wedge b].
\end{eqnarray*}
The proof is complete.
\end{proof}

   \begin{theorem}\label{th2}
 For all $M\in\mathcal {M}$ and   $\eta\in M_{G,A}^{2}(0,T)$, there exists a q.s. continuous
 version of stochastic integral
\begin{equation*}
\int_0^t \eta_{s}d M_{s}, \ \ 0\leq t\leq T,
\end{equation*}
i.e. there exists a continuous process $Y=\{Y_{t}\}_{t\in[0,T]}$ in
the sublinear expectation space $(\Omega, \mathcal {H},
\mathbb{\hat{E}})$ such that
\begin{equation*}
c(Y_{t}\neq\int_0^t \eta_{s}d M_{s})=0,\ \text{for all}\ t,  \ 0\leq
t\leq T.
\end{equation*}
\end{theorem}

\begin{proof}
We use $\pi^{n}=\{0=t_{0}^{n}<t_{1}^{n}\cdots<t_{n}^{n}=T\}$ to
denote a partition of $[0,T]$ such that
  $\max\{t_{i+1}^{n}-t_{i}^{n}, 0\leq i\leq n-1\}\rightarrow 0,$ as $n\rightarrow\infty.$

For any  $\eta\in M_{G,A}^{2}(0,T)$, there exists a sequence of
$\eta^{n}\in M_{G}^{2, 0}(0,T), n\geq 1,$ of the form
\begin{eqnarray*}
 \eta_{t}^{n}=\sum\limits_{j=0}^{n-1}\xi_{j}^{n}I_{[t_{j}^{n},t_{j+1}^{n})}(t),
 \end{eqnarray*}
 where $\xi_{j}^{n}\in L_{G}^{2}(\mathcal {F}_{t_{j}^{n}}),
 0\leq i\leq n-1$,
such that $$\mathbb{\hat{E}}[\int_0^T
 |\eta_{s}-\eta_{s}^{n}|^{2}dA_{s}]\rightarrow 0,\ \text{as}\ n\rightarrow\infty.$$
 We put
 $X_{t}^{n}=\int_0^t
 \eta_{s}^{n}d M_{s}=\sum\limits_{j=0}^{n-1}\xi_{j}^{n}( M_{t_{j+1}^{n}\wedge t}
 - M_{t_{j}^{n}\wedge t})$, for all $i=1,\cdots, n.$ Then $X^{n}$ is a
 continuous G-martingale and
 $$\mathbb{\hat{E}}[X_{t}^{n}|\mathcal {H}_{s}]=-\mathbb{\hat{E}}[-X_{t}^{n}|\mathcal
 {H}_{s}]=X_{s}^{n}, \ \text{for all}\ s\in[0,t].$$
 For any $\lambda>0$, by Markov inequality for capacity (see Lemma 13 in \cite{DHP:2008}) as well as
  Theorem \ref{th1} we have
\begin{eqnarray*}
&&c\Big(\{\sup\limits_{0\leq t\leq
T}|X_{t}^{n}-X_{t}^{m}|\geq\lambda\}\Big)
 \leq\frac{1}{\lambda^{2}}\mathbb{\hat{E}}[\sup\limits_{0\leq t\leq
T}|X_{t}^{n}-X_{t}^{m}|^{2}]\leq
\frac{4}{\lambda^{2}}\mathbb{\hat{E}}[|X_{T}^{n}-X_{T}^{m}|^{2}],
\end{eqnarray*}
and thanks to Proposition \ref{pr1} it follows that
\begin{eqnarray*}
&&c\Big(\{\sup\limits_{0\leq t\leq
T}|X_{t}^{n}-X_{t}^{m}|\geq\lambda\}\Big)\leq
 \frac{4}{\lambda^{2}}\mathbb{\hat{E}}\Big[\int_0^T
|\eta_{t}^{n}-\eta_{t}^{m}|^{2}dA_{t}\Big] \rightarrow\infty,
\end{eqnarray*}
as $n,m\rightarrow\infty$. Hence, we can choose a subsequence
$n_{k}\uparrow\infty$ such that

$$c\Big(\{\sup\limits_{0\leq t\leq
T}|X_{t}^{n_{k+1}}-X_{t}^{n_{k}}|\geq2^{-k}\}\Big)<2^{-k}, k\geq
1,$$ and from the Borel-Cantelli lemma for the capacity $c$ ( see
Lemma 5 in \cite{DHP:2008} ) we obtain
$$c\Big(\{\sup\limits_{0\leq t\leq
T}|X_{t}^{n_{k+1}}-X_{t}^{n_{k}}|\geq2^{-k},\ \text{for infinitely
many}\ k\}\Big)=0.$$
 Hence, there exists a random integer $k_{1}$ such that
$$\sup\limits_{0\leq t\leq
T}|X_{t}^{n_{k+1}}-X_{t}^{n_{k}}|<2^{-k},\ q.s., \ \text{for all}\
k>k_{1}.$$ This proves that the process $X$ converges  uniformly
  in $t\in[0,T]$ q.s.. We denote the limit of
$X^{n_{k}}$ by $Y$. Thanks to the quasi-sure uniform convergence it
is  a process whose paths are continuous. On the other hand, the
sequence $\{X_{t}^{n_{k}}, k\geq 1\}$ converges to $\int_0^t
\eta_{s}d M_{s}$ in $L^{2}_{G}(\mathcal {F})$, for all $t\in[0,T]$.
Thus,
\begin{equation*}
Y_{t}=\int_0^t \eta_{s}d M_{s},\ q.s.,  \ \text{for all}\ t,  \
0\leq t\leq T,
\end{equation*}
and so $Y$ is a continuous modification of the integral process. The
proof is complete.
\end{proof}

The following very useful Lemma was established by  Peng
\cite{Peng:2007}.
\begin{lemma}
Let $X,Y \in L^{1}_{G}(\mathcal {F}) $ be such that
$\mathbb{\hat{E}}[Y|\mathcal {H}_{s}]=-\mathbb{\hat{E}}[-Y|\mathcal
{H}_{s}]$, for $s\geq 0.$ Then we have
$$\mathbb{\hat{E}}[X+Y|\mathcal {H}_{s}]=\mathbb{\hat{E}}[X|\mathcal {H}_{s}]+\mathbb{\hat{E}}[Y|\mathcal
{H}_{s}].$$ In particular, if
$\mathbb{\hat{E}}[Y]=-\mathbb{\hat{E}}[-Y]=0$, then we have
$$\mathbb{\hat{E}}[X+Y]=\mathbb{\hat{E}}[X]+\mathbb{\hat{E}}[Y].$$
\end{lemma}

Now we give the Burkholder-Davis-Gundy inequality for the
stochastic integral with respect to G-martingales.
\begin{theorem}\label{th5}
For every $q>0$, there exist a positive constant $C_{q}$ such that,
for all $M\in\mathcal {M}$  and all $\eta\in M_{G,A}^{2}(0,T)$,
\begin{eqnarray*}
\mathbb{\hat{E}}\Big[\sup\limits_{t\in[0,T]}|\int_0^t \eta_{s}d
M_{s}|^{2q}\Big]\leq C_{q}\mathbb{\hat{E}}\Big[(\int_0^T
\eta_{s}^{2}dA_{s})^{q}\Big] .
\end{eqnarray*}
\end{theorem}

\begin{proof}
Let $M\in\mathcal {M}$. Then for all $P\in\mathcal {P}$, $M$ is a
continuous $P-$martingale and $M^{2}-A$ is a continuous
$P-$supermartingale.  Let $\langle M\rangle^{(P)}$ denote the
quadratic variation process of $M$ under $P$, i.e., the unique
continuous, $P-$predictable increasing process  $\langle
M\rangle^{(P)}$ such that $\langle M\rangle^{(P)}_{0}=0$ and
$M^{2}-\langle M\rangle^{(P)}$ is a $P-$martingale. Then
 $$\langle M\rangle^{(P)}-A=(M^{2}-A)-(M^{2}-\langle
 M\rangle^{(P)}), \langle M\rangle^{(P)}_{0}-A_{0}=0,$$
 is a continuous $P-$supermartingale and of finite variation.

Thanks to Doob-Meyer decomposition theorem, we have
$$\langle M\rangle^{(P)}-A=N^{(P)}-B^{(P)}, N^{(P)}_{0}-B^{(P)}_{0}=0,$$
where $N^{(P)}$ is a continuous $P-$martingale and $B^{(P)}$ is a
$P-$predictable, continuous  increasing process. Therefore, $\langle
M\rangle^{(P)}-A+B^{(P)}$  is a continuous $P-$martingale and of
finite variation. Consequently, $$\langle
M\rangle^{(P)}-A+B^{(P)}=0, P-a.s., i.e.,$$
\begin{eqnarray}\label{eq10}
 d\langle M\rangle^{(P)}_{t}\leq dA_{t}, t\geq 0, P-a.s..
\end{eqnarray}

Let $\eta\in M_{G,A}^{2}(0,T)$. Then for any $0\leq s \leq t\leq T$,
by the proof of Theorem \ref{th1},  we know that $E_{P}[\int_0^t
\eta_{r}d M_{r}|\mathcal {F}_{s}^{B}]= \int_0^s \eta_{r}d M_{r},
P$-$ a.s.,$  for all $P\in\mathcal {P}$.
 From the classical Burkh\"{o}lder-Davis-Gundy inequalities, for every $q>0$,
 there exist a positive constant $C_{q}$ such that
\begin{eqnarray*}
E_{P}\Big[\sup\limits_{t\in[0,T]}|\int_0^t \eta_{s}d
M_{s}|^{2q}\Big]\leq C_{q}E_{P}\Big[(\int_0^T \eta_{s}^{2}d\langle
M\rangle_{s}^{(P)})^{q}\Big],
\end{eqnarray*}
and from (\ref{eq10})
\begin{eqnarray*}
E_{P}\Big[\sup\limits_{t\in[0,T]}|\int_0^t \eta_{s}d
M_{s}|^{2q}\Big]\leq C_{q}E_{P}\Big[(\int_0^T \eta_{s}^{2}d
A_{s})^{q}\Big],
\end{eqnarray*}
Therefore,
\begin{eqnarray*}
\mathbb{\hat{E}}\Big[\sup\limits_{t\in[0,T]}|\int_0^t \eta_{s}d
M_{s}|^{2q}\Big]\leq C_{q}\mathbb{\hat{E}}\Big[(\int_0^T
\eta_{s}^{2}dA_{s})^{q}\Big] .
\end{eqnarray*}
 The proof is complete.
\end{proof}

Let $\pi^{n}=\{0=t_{0}^{n}<t_{1}^{n}\cdots<t_{n}^{n}=T\}$ with
$|\pi^{n}|\rightarrow 0,$ as $n\rightarrow\infty$, be a partition of
the interval $[0,T]$. In the following of this section, we assume
that the process $A$ satisfies the following assumption:

$\mathbb{\hat{E}}[A_{T}^{2}]<\infty,$ and for all
$\{\pi^{n}\}_{n\geq 1}$ sequence of partition of $[0,T]$ such that
$|\pi^{n}|\rightarrow 0,$ as $n\rightarrow\infty$,
$\mathbb{\hat{E}}[\sum\limits_{i=0}^{n-1}(A_{t_{i+1}^{n}}-A_{t_{i}^{n}})^{2}]\rightarrow
0, n\rightarrow\infty.$

\begin{proposition}\label{pr4}
Let $M\in\mathcal {M}$. Then the quadratic variation of $M$ exists
and
\begin{eqnarray*}
\langle M\rangle_{t}=M^{2}_{t}-2\int^{t}_{0}M_{s}dM_{s}, \ \text{for
all}\ t\geq 0.
\end{eqnarray*}
\end{proposition}

\begin{proof}
We use $\pi^{n}=\{0=t_{0}^{n}<t_{1}^{n}\cdots<t_{n}^{n}=T\}$ to
denote a partition of $[0,T]$ such that
  $\max\{t_{i+1}^{n}-t_{i}^{n}, 0\leq i\leq n-1\}\rightarrow 0,$ as
  $n\rightarrow\infty.$ Then
\begin{eqnarray}\label{eq6}
M^{2}_{t}&=&
\sum\limits_{i=0}^{n-1}[M_{t_{i+1}^{n}\wedge t}^{2}-M_{t_{i}^{n}\wedge t}^{2}]\nonumber\\
&=&2\sum\limits_{i=0}^{n-1}M_{t_{i}^{n}\wedge
t}[M_{t_{i+1}^{n}\wedge t}-M_{t_{i}^{n}\wedge t}]
+\sum\limits_{i=0}^{n-1}[M_{t_{i+1}^{n}\wedge t}-M_{t_{i}^{n}\wedge
t}]^{2}.
\end{eqnarray}
Thanks to Theorem \ref{th5}, we have
\begin{eqnarray*}
&&\mathbb{\hat{E}}[\int^{T}_{0}(M_{s}-M_{s}^{n})^{2}dA_{s}]\\
&=&C\mathbb{\hat{E}}[\sum\limits_{i=0}^{n-1}\int^{t_{i+1}^{n}}_{t_{i}^{n}}(M_{s}-M_{s}^{n})^{2}dA_{s}]\\
&\leq&C\mathbb{\hat{E}}[\sum\limits_{i=0}^{n-1}\sup\limits_{s\in[t_{i}^{n},t_{i+1}^{n}]}
(M_{s}-M_{t_{i}^{n}})^{2}(A_{t_{i+1}^{n}}-A_{t_{i}^{n}})]\\
&\leq&C\mathbb{\hat{E}}[\sum\limits_{i=0}^{n-1}\sup\limits_{t\in[0,T]}|\int^{t}_{0}
I_{[t_{i}^{n},t_{i+1}^{n}]}(s)dM_{s}|^{4}]^{\frac{1}{2}}
\mathbb{\hat{E}}[\sum\limits_{i=0}^{n-1}(A_{t_{i+1}^{n}}-A_{t_{i}^{n}})^{2}]^{\frac{1}{2}}\\
&\leq&C\mathbb{\hat{E}}[\sum\limits_{i=0}^{n-1}(A_{t_{i+1}^{n}}-A_{t_{i}^{n}})^{2}]\rightarrow
0, \text{as}\ n\rightarrow \infty.
\end{eqnarray*}
Therefore,
\begin{eqnarray}\label{e1}
&&\mathbb{\hat{E}}[\sup\limits_{t\in[0,T]}|\int^{t}_{0}(M_{s}-M_{s}^{n})dM_{s}|^{2}]\leq
C\mathbb{\hat{E}}[\int^{T}_{0}(M_{s}-M_{s}^{n})^{2}dA_{s}]\rightarrow
0, \text{as}\ n\rightarrow \infty.
\end{eqnarray}
Consequently,  the first term of (\ref{eq6}) converges to the
stochastic integral $2\int^{t}_{0}M_{s}dM_{s}$, then the quadratic
variation of $M$ exists and is equal to
\begin{eqnarray*}
\langle
M\rangle_{t}:=\lim\limits_{n\rightarrow\infty}\sum\limits_{i=0}^{n-1}[M_{t_{i+1}^{n}\wedge
t}-M_{t_{i}^{n}\wedge t}]^{2} =M^{2}_{t}-2\int^{t}_{0}M_{s}dM_{s}.
\end{eqnarray*}
The proof is complete.
\end{proof}

By Theorem \ref{th5} and Proposition \ref{pr4}, we have
\begin{remark}\label{re1}
For all $t\in [0,T]$, we have $\langle M\rangle_{t}^{(P)}=\langle
M\rangle_{t}, P-a.s.,$ for all $P\in\mathcal {P}$.
\end{remark}

\begin{definition}
Let $M\in\mathcal {M}$. Then for all $\eta\in M_{G}^{1,0}(0,T)$ of
the form
 $\eta_{t}=\sum\limits_{j=0}^{n-1}\xi_{t_{j}}I_{[t_{j},t_{j+1})}(t)$
  we define $$I(\eta)=\int_0^T
 \eta_{t}d\langle M\rangle_{t}=\sum\limits_{j=0}^{n-1}\xi_{t_{j}}(\langle M\rangle_{t_{j+1}}
 -\langle M\rangle_{t_{j}}).$$
\end{definition}

We have the following proposition.
\begin{proposition}\label{pr2}
For any $M\in\mathcal {M}$, the mapping $I:
M_{G}^{1,0}(0,T)\rightarrow L_{G}^{1}(\mathcal {F}_{T})$ is a linear
continuous mapping, and thus, can be continuously extended to $I:
M_{G,A}^{1}(0,T)\rightarrow L_{G}^{1}(\mathcal {F}_{T})$. Moreover,
for all $\eta\in M_{G,A}^{1}(0,T)$ we have
\begin{eqnarray}\label{ineq2}
\mathbb{\hat{E}}\Big[|\int_0^T \eta_{t}d\langle M\rangle_{t}|\Big]
 \leq \mathbb{\hat{E}}\Big[\int_0^T |\eta_{t}|dA_{t}\Big].
\end{eqnarray}
\end{proposition}
\begin{proof}
From $M$ is a symmetric G-martingale and $M^{2}-A\ $ is a
G-supermartingale it follows that, for all $0\leq s \leq t \leq T$,
\begin{eqnarray*}
&&\mathbb{\hat{E}}[\langle M\rangle_{t}-\langle M\rangle_{s}-(A_{t}-A_{s})|\mathcal {H}_{s}]\\
&=&\mathbb{\hat{E}}[M_{t}^{2}-M_{s}^{2}-2\int^{t}_{s}M_{r}dM_{r}-(A_{t}-A_{s})|\mathcal
{H}_{s}]\\
&=&\mathbb{\hat{E}}[M_{t}^{2}-M_{s}^{2}-(A_{t}-A_{s})|\mathcal
{H}_{s}]\leq 0.
\end{eqnarray*}
For $\eta\in M_{G}^{1,0}(0,T)$ of the form
 $\eta_{t}=\sum\limits_{j=0}^{n-1}\xi_{t_{j}}I_{[t_{j},t_{j+1})}(t),$ we
 have
\begin{eqnarray*}
 &&\mathbb{\hat{E}}\Big[|\int_0^T \eta_{t}d\langle M\rangle_{t}|\Big]
 = \mathbb{\hat{E}}\Big[|\sum\limits_{j=0}^{N-1}\xi_{t_{j}}(\langle M\rangle_{t_{j+1}}
 -\langle M\rangle_{t_{j}})|\Big]\\
&\leq&
\mathbb{\hat{E}}\Big[\sum\limits_{j=0}^{N-1}|\xi_{t_{j}}|(\langle
M\rangle_{t_{j+1}}-\langle M\rangle_{t_{j}})\Big]\\
   &\leq& \mathbb{\hat{E}}\Big[\sum\limits_{j=0}^{N-1}|\xi_{t_{j}}|(\langle
M\rangle_{t_{j+1}}-\langle M\rangle_{t_{j}}-A_{t_{j+1}}
 +A_{t_{j}})\Big]+\mathbb{\hat{E}}\Big[\sum\limits_{j=0}^{N-1}|\xi_{t_{j}}|(A_{t_{j+1}}
 -A_{t_{j}})\Big],
 \end{eqnarray*}
  where

 \begin{eqnarray*}
&& \mathbb{\hat{E}}\Big[\sum\limits_{j=0}^{N-1}|\xi_{t_{j}}|(\langle
M\rangle_{t_{j+1}}-\langle M\rangle_{t_{j}}-A_{t_{j+1}}
 +A_{t_{j}})\Big]\\
  &\leq& \mathbb{\hat{E}}\Big[\sum\limits_{j=0}^{N-2}|\xi_{t_{j}}|(\langle
M\rangle_{t_{j+1}}-\langle M\rangle_{t_{j}}-A_{t_{j+1}}
 +A_{t_{j}}) + \\&&|\xi_{t_{N-1}}|  \mathbb{\hat{E}}[\langle
M\rangle_{t_{N}}-\langle M\rangle_{t_{N-1}}-A_{t_{N}}
 +A_{t_{N-1}}|\mathcal {H}_{t_{N-1}}]\Big]\\
 &\leq& \mathbb{\hat{E}}\Big[\sum\limits_{j=0}^{N-2}|\xi_{t_{j}}|(\langle
M\rangle_{t_{j+1}}-\langle M\rangle_{t_{j}}-A_{t_{j+1}}
 +A_{t_{j}})\Big] \leq \cdots\leq 0.
 \end{eqnarray*}
 Consequently,
  \begin{eqnarray*}
 \mathbb{\hat{E}}\Big[|\int_0^T \eta_{t}d\langle M\rangle_{t}|\Big]
 &\leq& \mathbb{\hat{E}}\Big[\sum\limits_{j=0}^{N-1}|\xi_{t_{j}}|(A_{t_{j+1}}
 -A_{t_{j}})\Big]=\mathbb{\hat{E}}\Big[\int_0^T
 |\eta_{t}|dA_{t}\Big].
 \end{eqnarray*}
 Thus, (\ref{ineq2}) holds for all $\eta\in M_{G}^{1,0}(0,T)$. We
then can continuously extend the above inequality to the case
$\eta\in M_{G,A}^{1}(0,T)$ and prove (\ref{ineq2}). The proof is
complete.
\end{proof}

Now we can  prove the following  proposition by the same argument as
in \cite{Peng:2007}.

\begin{proposition}\label{pr2}
If  $M\in\mathcal {M}$,  $X\in L^{1}_{G}(\mathcal {F})$ and $\xi\in
L^{2}_{G}(\mathcal {F}_{s})$, then for all $0\leq s \leq t <\infty$
\begin{eqnarray*}
\mathbb{\hat{E}}[X+\xi(\langle M\rangle_{t}-\langle
M\rangle_{s})]=\mathbb{\hat{E}}[X+\xi( M_{t}-
M_{s})^{2}]=\mathbb{\hat{E}}[X+\xi( M_{t}^{2}- M_{s}^{2})].
\end{eqnarray*}
\end{proposition}

Moreover, we have the following isometry property.
 \begin{proposition}\label{pr3}
 If $M\in\mathcal {M}$  and  $\eta\in M_{G,A}^{2}(0,T)$, then
\begin{eqnarray*}
\mathbb{\hat{E}}\Big[(\int_0^T \eta_{t}d
M_{t})^{2}\Big]=\mathbb{\hat{E}}\Big[\int_0^T \eta_{t}^{2}d\langle
M\rangle_{t}\Big] .
\end{eqnarray*}
\end{proposition}

\begin{proof}
For $\eta\in M_{G}^{2,0}(0,T)$ of the form
 $\eta_{t}=\sum\limits_{j=0}^{n-1}\xi_{t_{j}}I_{[t_{j},t_{j+1})}(t)$
 a straightforward argument gives
\begin{eqnarray*}
 \mathbb{\hat{E}}\Big[|\int_0^T \eta_{t}dM_{t}|^{2}\Big]
 &= &\mathbb{\hat{E}}\Big[|\sum\limits_{j=0}^{n-1}\xi_{t_{j}}(M_{t_{j+1}}
 -M_{t_{j}})|^{2}\Big]
=\mathbb{\hat{E}}\Big[\sum\limits_{j=0}^{n-1}\xi_{t_{j}}^{2}(M_{t_{j+1}}
 -M_{t_{j}})^{2}\Big].
  \end{eqnarray*}
  Thanks to Proposition \ref{pr2}, we have
\begin{eqnarray*}
 \mathbb{\hat{E}}\Big[|\int_0^T \eta_{t}dM_{t}|^{2}\Big]
 =\mathbb{\hat{E}}\Big[\sum\limits_{j=0}^{n-1}\xi_{t_{j}}^{2}(\langle
M\rangle_{t_{j+1}}-\langle
M\rangle_{t_{j}})\Big]=\mathbb{\hat{E}}\Big[\int_0^T
\eta_{t}^{2}d\langle M\rangle_{t}\Big].
\end{eqnarray*}
 Thus, Proposition \ref{pr3} holds for all $\eta\in M_{G}^{2,0}(0,T)$. Finally,
 Proposition \ref{pr2} allows to extend continuously  the above inequality to
 all $\eta\in M_{G,A}^{2}(0,T)$, and thus, yields the desired result. The proof is
complete.
\end{proof}

Now we give another kind of the Burkholder-Davis-Gundy inequalities
for the stochastic integral with respect to G-martingales.
\begin{theorem}\label{th12}
For every $p>0$, there exist two positive constants $c_{p}$ and
$C_{p}$ such that, for all $M\in\mathcal {M}$ such that $M\in
M_{G,A}^{2}(0,T)$ and all $\eta\in M_{G,A}^{2}(0,T)$,
\begin{eqnarray*}
c_{p}\mathbb{\hat{E}}\Big[(\int_0^T \eta_{s}^{2}d\langle
M\rangle_{s})^{p}\Big]\leq
\mathbb{\hat{E}}\Big[\sup\limits_{t\in[0,T]}|\int_0^t \eta_{s}d
M_{s}|^{2p}\Big]\leq C_{p}\mathbb{\hat{E}}\Big[(\int_0^T
\eta_{s}^{2}d\langle M\rangle_{s})^{p}\Big] .
\end{eqnarray*}
\end{theorem}

\begin{proof}
 For any $0\leq s
\leq t\leq T$, by the proof of Theorem \ref{th1},  we know that
$\int_0^\cdot \eta_{r}d M_{r}$ is a continuous $P-$martingale,  for
all $P\in\mathcal {P}$. Thanks to Remark \ref{re1}, we know that
$\langle M\rangle_{t}^{(P)}=\langle M\rangle_{t}, P-a.s.$ all $t\in
[0,T]$.
 From the classical Burkh\"{o}lder-Davis-Gundy inequalities it follows that
\begin{eqnarray*}
c_{p}E_{P}\Big[(\int_0^T \eta_{s}^{2}d\langle
M\rangle_{s})^{p}\Big]\leq
E_{P}\Big[\sup\limits_{t\in[0,T]}|\int_0^t \eta_{s}d
M_{s}|^{2p}\Big]\leq C_{p}E_{P}\Big[(\int_0^T \eta_{s}^{2}d\langle
M\rangle_{s})^{p}\Big],
\end{eqnarray*}
where two constants $0\leq c_{p} \leq C_{p}$ only depend on $p$.
Therefore,
\begin{eqnarray*}
c_{p}\mathbb{\hat{E}}\Big[(\int_0^T \eta_{s}^{2}d\langle
M\rangle_{s})^{p}\Big]\leq
\mathbb{\hat{E}}\Big[\sup\limits_{t\in[0,T]}|\int_0^t \eta_{s}d
M_{s}|^{2p}\Big]\leq C_{p}\mathbb{\hat{E}}\Big[(\int_0^T
\eta_{s}^{2}d\langle M\rangle_{s})^{p}\Big].
\end{eqnarray*}
 The proof is complete.
\end{proof}

 \begin{theorem}\label{th3}
 For all $t\geq 0$. Let $M\in\mathcal {M}$.
 If  $f\in M_{G,A}^{2}(0,t)$ is a bounded process, then
 the quadratic variation process of $X_{t}:=\int_0^t f_{s}d
M_{s}$ exists and
\begin{eqnarray*}
\langle X\rangle_{t}=\int_0^t f_{s}^{2}d\langle M\rangle_{s}.
\end{eqnarray*}
\end{theorem}

\begin{proof}
For $f\in M_{G,A}^{2}(0,t)$, there exists an $f^{n}$ of the form $
f_{s}^{n}=\sum\limits_{j=0}^{n-1}\xi_{t_{j}^{n}}I_{[t_{j}^{n},t_{j+1}^{n})}(s),$
 where $\xi_{t_{j}^{n}}\in L_{G}^{2}(\mathcal {F}_{t_{j}^{n}}),
 0\leq i\leq n-1$,
such that
 \begin{eqnarray*}
 \mathbb{\hat{E}}\Big[(\int_0^t (f_{s}-f_{s}^{n})d
M_{s})^{2}\Big]\leq \mathbb{\hat{E}}[\int_0^t
|f_{s}-f_{s}^{n}|^{2}dA_{s}]\rightarrow 0,\ \text{as}\
n\rightarrow\infty.
\end{eqnarray*}

For any $\varepsilon>0$ and $0\leq s \leq t$, we have
\begin{eqnarray*}
 \mathbb{\hat{E}}[|\int_s^t f_{r}^{2}dA_{r}-\int_s^t
(f_{r}^{n})^{2}dA_{r}|]&\leq& (1+\frac{1}{\varepsilon})
\mathbb{\hat{E}}[\int_s^t |f_{r}-f_{r}^{n}|^{2}dA_{r}]+ \varepsilon
\mathbb{\hat{E}}[\int_s^t f_{r}^{2}dA_{r}]\\
& \rightarrow&  \varepsilon \mathbb{\hat{E}}[\int_s^t
f_{r}^{2}dA_{r}], \ \text{as}\ n\rightarrow\infty.
\end{eqnarray*}
Therefore,
\begin{eqnarray*}
 \mathbb{\hat{E}}[|\int_s^t f_{r}^{2}dA_{r}-\int_s^t
(f_{r}^{n})^{2}dA_{r}|] \rightarrow 0, \ \text{as}\
n\rightarrow\infty.
\end{eqnarray*}
Since $M^{2}-A$ is a G-supermartingale, we have
\begin{eqnarray*}
&& \mathbb{\hat{E}}\Big[(\int_s^t f_{r}^{n}d M_{r})^{2}-\int_s^t
(f_{r}^{n})^{2}dA_{r}|\mathcal {H}_{s}\Big]\\
&=&\mathbb{\hat{E}}[(\sum\limits_{j=0}^{n-1}\xi_{t_{j}^{n}}(M_{t_{j+1}^{n}\vee
s}-M_{t_{j}^{n}\vee s}))^{2}-
\sum\limits_{j=0}^{n-1}\xi^{2}_{t_{j}^{n}}(A_{t_{j+1}^{n}\vee
s}-A_{t_{j}^{n}\vee s})|\mathcal
{H}_{s}]\\
&=&\mathbb{\hat{E}}[\sum\limits_{j=0}^{n-1}\xi^{2}_{t_{j}^{n}}(M_{t_{j+1}^{n}\vee
s}-M_{t_{j}^{n}\vee s})^{2}-
\sum\limits_{j=0}^{n-1}\xi^{2}_{t_{j}^{n}}(A_{t_{j+1}^{n}\vee
s}-A_{t_{j}^{n}\vee s})|\mathcal
{H}_{s}]\\
&=&\mathbb{\hat{E}}[\sum\limits_{j=0}^{n-1}\xi^{2}_{t_{j}^{n}}(M_{t_{j+1}^{n}\vee
s}^{2}-M_{t_{j}^{n}\vee s}^{2}- A_{t_{j+1}^{n}\vee
s}+A_{t_{j}^{n}\vee s})|\mathcal {H}_{s}]\\
&\leq&\sum\limits_{j=0}^{n-1}\mathbb{\hat{E}}[\xi^{2}_{t_{j}^{n}}(M_{t_{j+1}^{n}\vee
s}^{2}-M_{t_{j}^{n}\vee s}^{2}- A_{t_{j+1}^{n}\vee
s}+A_{t_{j}^{n}\vee s})|\mathcal {H}_{s}]\leq 0.
\end{eqnarray*}
For all $\varepsilon>0$, from the above inequalities it follows that
\begin{eqnarray*}
 &&\mathbb{\hat{E}}\Big[\Big(\mathbb{\hat{E}}[(\int_s^t f_{r}d M_{r})^{2}-\int_s^t
f_{r}^{2}dA_{r}|\mathcal {H}_{s}]\Big)^{+}\Big]\\&\leq&
\mathbb{\hat{E}}\Big[\Big(\mathbb{\hat{E}}[(\int_s^t f_{r}d
M_{r})^{2}-\int_s^t f_{r}^{2}dA_{r}|\mathcal {H}_{s}]-
\mathbb{\hat{E}}[(\int_s^t f_{r}^{n}d M_{r})^{2}-\int_s^t
(f_{r}^{n})^{2}dA_{r}|\mathcal {H}_{s}]\Big)^{+}\Big]
\\&&+\mathbb{\hat{E}}\Big[\Big(
\mathbb{\hat{E}}\Big[(\int_s^t f_{r}^{n}d M_{r})^{2}-\int_s^t
(f_{r}^{n})^{2}dA_{r}|\mathcal {H}_{s}\Big]\Big)^{+}\Big]
\\&=&
\mathbb{\hat{E}}\Big[\Big(\mathbb{\hat{E}}[(\int_s^t f_{r}d
M_{r})^{2}-\int_s^t f_{r}^{2}dA_{r}|\mathcal {H}_{s}]-
\mathbb{\hat{E}}[(\int_s^t f_{r}^{n}d M_{r})^{2}-\int_s^t
(f_{r}^{n})^{2}dA_{r}|\mathcal {H}_{s}]\Big)^{+}\Big]
\\&\leq&
\mathbb{\hat{E}}\Big[\Big|\mathbb{\hat{E}}[(\int_s^t f_{r}d
M_{r})^{2}-\int_s^t f_{r}^{2}dA_{r}|\mathcal {H}_{s}]-
\mathbb{\hat{E}}[(\int_s^t f_{r}^{n}d M_{r})^{2}-\int_s^t
(f_{r}^{n})^{2}dA_{r}|\mathcal {H}_{s}]\Big|\Big]
\\&\leq&
\mathbb{\hat{E}}\Big[\Big|(\int_s^t f_{r}d M_{r})^{2}- (\int_s^t
f_{r}^{n}d M_{r})^{2}\Big|\Big]+\mathbb{\hat{E}}\Big[\Big|\int_s^t
f_{r}^{2}dA_{r}-\int_s^t (f_{r}^{n})^{2}dA_{r}\Big|\Big]
\\&\leq&
(1+\frac{1}{\varepsilon}) \mathbb{\hat{E}}\Big[(\int_s^t
(f_{r}-f_{r}^{n})d M_{r})^{2}\Big]+
\varepsilon\mathbb{\hat{E}}\Big[(\int_s^t f_{r}d M_{r})^{2}\Big]+
\mathbb{\hat{E}}\Big[\Big|\int_s^t f_{r}^{2}dA_{r}-\int_s^t
(f_{r}^{n})^{2}dA_{r}\Big|\Big]\\
&\rightarrow& \varepsilon\mathbb{\hat{E}}\Big[(\int_s^t f_{r}d
M_{r})^{2}\Big],  \ \text{as}\ n\rightarrow\infty.
\end{eqnarray*}
Therefore, $$ \mathbb{\hat{E}}\Big[\Big(\mathbb{\hat{E}}[(\int_s^t
f_{r}d M_{r})^{2}-\int_s^t f_{r}^{2}dA_{r}|\mathcal
{H}_{s}]\Big)^{+}\Big]=0,$$ and which yields
$$\mathbb{\hat{E}}[(\int_s^t
f_{r}d M_{r})^{2}-\int_s^t f_{r}^{2}dA_{r}|\mathcal {H}_{s}]\leq0, \
q.s.$$ From the above inequality and Proposition \ref{pr2} it
follows that for all $0\leq s \leq t$
\begin{eqnarray*}
&&\mathbb{\hat{E}}[X_{t}^{2}-\int_0^t f_{r}^{2}dA_{r}|\mathcal
{H}_{s}]=\mathbb{\hat{E}}[(X_{s}+\int_s^t f_{r}d M_{r})^{2}-\int_0^t
f_{r}^{2}dA_{r}|\mathcal{H}_{s}]
\\&=&\mathbb{\hat{E}}[X_{s}^{2}+2X_{s}\int_s^t f_{r}d
M_{r}+(\int_s^t f_{r}d M_{r})^{2}-\int_0^t f_{r}^{2}dA_{r}|\mathcal
{H}_{s}]
\\&=&X_{s}^{2}-\int_0^s f_{r}^{2}dA_{r}
+\mathbb{\hat{E}}[(\int_s^t f_{r}d M_{r})^{2}-\int_s^t
f_{r}^{2}dA_{r}|\mathcal {H}_{s}] \leq X_{s}^{2}-\int_0^s
f_{r}^{2}dA_{r}.
\end{eqnarray*}
Consequently,  $X_{t}^{2}-\int_0^t f_{r}^{2}dA_{r}$ is a
G-supermartingale.

From Proposition \ref{pr1}, we know that $X$ is a symmetric
 G-martingale. Then from Proposition \ref{pr3} it follows that the quadratic variation process of $X$
 exists.
Therefore
\begin{eqnarray*}
\mathbb{\hat{E}}[|\langle X\rangle_{t}-\int_0^t f_{s}^{2}d\langle
M\rangle_{s}|]&\leq&  \mathbb{\hat{E}}[|\langle
X\rangle_{t}-\sum\limits_{i=0}^{n-1}(X_{t_{i+1}^{n}}-X_{t_{i}^{n}})^{2}|]\\
&&+\mathbb{\hat{E}}[|\sum\limits_{i=0}^{n-1}(X_{t_{i+1}^{n}}-X_{t_{i}^{n}})^{2}
-\sum\limits_{i=0}^{n-1}\xi_{t_{i}^{n}}^{2}(M_{t_{i+1}^{n}}-M_{t_{i}^{n}})^{2}|]\\
&&+\mathbb{\hat{E}}[|\sum\limits_{i=0}^{n-1}\xi_{t_{i}^{n}}^{2}(M_{t_{i+1}^{n}}-M_{t_{i}^{n}})^{2}-
\sum\limits_{i=0}^{n-1}\xi_{t_{i}^{n}}^{2}(\langle M\rangle _{t_{i+1}^{n}}-\langle M\rangle _{t_{i}^{n}})|]\\
&&+\mathbb{\hat{E}}[|\sum\limits_{i=0}^{n-1}\xi_{t_{i}^{n}}^{2}(\langle
M\rangle _{t_{i+1}^{n}}-\langle M\rangle _{t_{i}^{n}})-\int_0^t
f_{s}^{2}d\langle M\rangle_{s}|]\\
&:=&I_{1}+I_{2}+I_{3}+I_{4}.
\end{eqnarray*}

As $n\rightarrow \infty$, $I_{1}\rightarrow 0, I_{4}\rightarrow 0$.
Now we prove  $I_{2}\rightarrow 0, I_{3}\rightarrow 0$, as
$n\rightarrow \infty$.

\begin{eqnarray*}
I_{2}&=&\mathbb{\hat{E}}[|\sum\limits_{i=0}^{n-1}(\int_{t_{i}^{n}}^{t_{i+1}^{n}}
f_{s}dM_{s})^{2}-\sum\limits_{i=0}^{n-1}(\int_{t_{i}^{n}}^{t_{i+1}^{n}}
f_{s}^{n}dM_{s})^{2}|]\\
&\leq&\mathbb{\hat{E}}[\sum\limits_{i=0}^{n-1}|(\int_{t_{i}^{n}}^{t_{i+1}^{n}}
(f_{s}-f_{s}^{n})dM_{s})^{2}-2(\int_{t_{i}^{n}}^{t_{i+1}^{n}}
f_{s}dM_{s})(\int_{t_{i}^{n}}^{t_{i+1}^{n}}
(f_{s}-f_{s}^{n})dM_{s})|]\\
&\leq&\mathbb{\hat{E}}[(1+\frac{1}{\varepsilon})\sum\limits_{i=0}^{n-1}(\int_{t_{i}^{n}}^{t_{i+1}^{n}}
(f_{s}-f_{s}^{n})dM_{s})^{2}+\varepsilon\sum\limits_{i=0}^{n-1}
(\int_{t_{i}^{n}}^{t_{i+1}^{n}}
f_{s}dM_{s})^{2}]\\
&\leq&(1+\frac{1}{\varepsilon})\mathbb{\hat{E}}[\sum\limits_{i=0}^{n-1}(\int_{t_{i}^{n}}^{t_{i+1}^{n}}
(f_{s}-f_{s}^{n})dM_{s})^{2}]+\varepsilon\mathbb{\hat{E}}[\sum\limits_{i=0}^{n-1}(\int_{t_{i}^{n}}^{t_{i+1}^{n}}
f_{s}dM_{s})^{2}]\\
&=&(1+\frac{1}{\varepsilon})\mathbb{\hat{E}}[(\int_0^t
(f_{s}-f_{s}^{n})dM_{s})^{2}]+\varepsilon\mathbb{\hat{E}}[(\int_0^t
f_{s}dM_{s})^{2}]
\end{eqnarray*}
Thanks to Proposition \ref{pr1}, we have
\begin{eqnarray*}
I_{2} &\leq&(1+\frac{1}{\varepsilon})\mathbb{\hat{E}}[\int_{0}^{t}
(f_{s}-f_{s}^{n}))^{2}dA_{s}]+\varepsilon\mathbb{\hat{E}}[\int_{0}^{t}
f_{s}^{2}dA_{s}]\\
&\rightarrow&\varepsilon\mathbb{\hat{E}}[\int_{0}^{t}
f_{s}^{2}dA_{s}],\ \text{as}\  n\rightarrow\infty,
\end{eqnarray*}
where
 $M_{s}^{n}=\sum\limits_{j=0}^{n-1}M_{t_{j}^{n}}I_{[t_{j}^{n},t_{j+1}^{n})}(s).$

From Proposition \ref{pr1} and the properties of G-expectation it
follows that
\begin{eqnarray*}
&&\mathbb{\hat{E}}[|\sum\limits_{i=0}^{n-1}\xi_{t_{i}^{n}}^{2}(M_{t_{i+1}^{n}}-M_{t_{i}^{n}})^{2}-
\sum\limits_{i=0}^{n-1}\xi_{t_{i}^{n}}^{2}(\langle M\rangle _{t_{i+1}^{n}}-\langle M\rangle _{t_{i}^{n}})|^{2}]\\
&=&4\mathbb{\hat{E}}[|\sum\limits_{i=0}^{n-1}\xi_{t_{i}^{n}}^{2}\int_{t_{i}^{n}}^{t_{i+1}^{n}}
(M_{s}-M_{t_{i}^{n}})d M_{s}|^{2}]\\
&=&4\mathbb{\hat{E}}[\sum\limits_{i=0}^{n-1}\xi_{t_{i}^{n}}^{4}(\int_{t_{i}^{n}}^{t_{i+1}^{n}}
(M_{s}-M_{t_{i}^{n}})d M_{s})^{2}]\\
&\leq&4C\mathbb{\hat{E}}[\sum\limits_{i=0}^{n-1}(\int_{t_{i}^{n}}^{t_{i+1}^{n}}
(M_{s}-M_{t_{i}^{n}})d M_{s})^{2}]\\
&=&4C\mathbb{\hat{E}}[(\int_0^t (M_{s}-M_s^{n})d
M_{s})^{2}]\rightarrow 0, \ \text{as}\ n \rightarrow\infty.
\end{eqnarray*}
Therefore $I_{3}\rightarrow 0, \ \text{as}\ n \rightarrow\infty.$

By the above inequality, we have
\begin{eqnarray*}
\mathbb{\hat{E}}[|\langle X\rangle_{t}-\int_0^t f_{s}^{2}d\langle
M\rangle_{s}|]\leq \varepsilon\int_{0}^{t}
\mathbb{\hat{E}}[f_{s}^{2}]dA_{s}.
\end{eqnarray*}
Thus,
\begin{eqnarray*}
\langle X\rangle_{t}=\int_0^t f_{s}^{2}d\langle M\rangle_{s}, \ q.s.
\end{eqnarray*}
We obtain the desired result. The proof is completed.
\end{proof}

\begin{proposition}\label{}
 For a fixed $T\geq 0$, $M$ is a symmetric G-martingale, $M^{2}-A\ $ is a
G-martingale and for $0\leq\sigma^{2}_{0}\leq 1,$
$-(M^{2}-\sigma^{2}_{0}A)$ is a G-martingale, if  $f\in
M_{G,A}^{1}(0,T)$,
  then
 the  process
\begin{eqnarray*}
X_{t}:=\int_0^t f_{s}d\langle M\rangle_{s}-2\int_0^t G(f_{s})d
A_{s}, \ t\in[0,T]
\end{eqnarray*}
is a decreasing G-martingale.
\end{proposition}

\begin{proof}
It is easy to check that $X$ is a decreasing G-martingale.
   We use $\pi=\{0=t_{0}^{n}<t_{1}^{n}\cdots<t_{n}^{n}=T\}$ to
denote a partition of $[0,T]$ such that
  $\max\{t_{i+1}^{n}-t_{i}^{n}, 0\leq i\leq n-1\}\rightarrow 0,$ as $n\rightarrow\infty.$

For $f\in M_{G,A}^{1}(0,T)$, there exists an $f^{n}$ of the form $
f_{t}^{n}=\sum\limits_{j=0}^{n-1}\xi_{t_{j}}I_{[t_{j}^{n},t_{j+1}^{n})}(t),$
 where $\xi_{t_{j}}\in L_{G}^{1}(\mathcal {F}_{t_{j}^{n}}),
 0\leq i\leq n-1$.
Let
\begin{eqnarray*}
X_{t}^{n}:=\sum\limits_{i=0}^{n-1}\xi_{t_{i}^{n}}(\langle M\rangle
_{t_{i+1}^{n}\wedge t}-\langle M\rangle _{t_{i}^{n}\wedge
t})-2\sum\limits_{i=0}^{n-1}G(\xi_{t_{i}^{n}})( A
_{t_{i+1}^{n}\wedge t}-A _{t_{i}^{n}\wedge t}),
\end{eqnarray*}
where $t\in[0,T]$.

 For $0\leq s \leq t \leq T$. Without loss of generality, we
suppose that $t_{k-1}^{n}\leq s \leq t_{k}^{n}  \leq t \leq
t_{k+1}^{n}$, for some $k=1,\cdots, n-1.$ Thus,
\begin{eqnarray*}
\mathbb{\hat{E}}[X_{t}^{n}|\mathcal
{H}_{t_{k}^{n}}]&=&\mathbb{\hat{E}}[\sum\limits_{i=0}^{n-1}\xi_{t_{i}^{n}}(\langle
M\rangle _{t_{i+1}^{n}\wedge t}-\langle M\rangle _{t_{i}^{n}\wedge
t})-2\sum\limits_{i=0}^{n-1}G(\xi_{t_{i}^{n}})( A
_{t_{i+1}^{n}\wedge t}-A _{t_{i}^{n}\wedge t})|\mathcal
{H}_{t_{k}^{n}}]\\
&=&\mathbb{\hat{E}}[\sum\limits_{i=0}^{k-1}\xi_{t_{i}^{n}}(\langle
M\rangle _{t_{i+1}^{n}}-\langle M\rangle
_{t_{i}^{n}})-2\sum\limits_{i=0}^{k-1}G(\xi_{t_{i}^{n}})( A
_{t_{i+1}^{n}}-A _{t_{i}^{n}})\\&&+\xi_{t_{k}^{n}}(\langle M\rangle
_{t}-\langle M\rangle _{t_{k}^{n}})-2G(\xi_{t_{k}^{n}})( A _{t}-A
_{t_{k}^{n}})|\mathcal
{H}_{t_{k}^{n}}]\\
&=&X_{t_{k}^{n}}^{n}+\mathbb{\hat{E}}[\xi_{t_{k}^{n}}(\langle
M\rangle _{t}-\langle M\rangle _{t_{k}^{n}})-2G(\xi_{t_{k}^{n}})(A
_{t}-A _{t_{k}^{n}})|\mathcal
{H}_{t_{k}^{n}}]\\
&=&X_{t_{k}^{n}}^{n}+\xi_{t_{k}^{n}}^{+}\mathbb{\hat{E}}[\langle
M\rangle _{t}-\langle M\rangle _{t_{k}^{n}}-A _{t}+A
_{t_{k}^{n}}|\mathcal
{H}_{t_{k}^{n}}]\\&&+\xi_{t_{k}^{n}}^{-}\mathbb{\hat{E}}[-\langle
M\rangle _{t}+\langle M\rangle _{t_{k}^{n}}+\sigma_{0}^{2}( A _{t}-A
_{t_{k}^{n}})|\mathcal
{H}_{t_{k}^{n}}]\\
&=&X_{t_{k}^{n}}^{n}.
\end{eqnarray*}

From Proposition \ref{p2} it follows that
\begin{eqnarray*}
\mathbb{\hat{E}}[X_{t}^{n}|\mathcal
{H}_{s}]&=&\mathbb{\hat{E}}[\mathbb{\hat{E}}[X_{t}^{n}|\mathcal
{H}_{t_{k}^{n}}]|\mathcal{H}_{s}]=\mathbb{\hat{E}}[X_{t_{k}^{n}}^{n}|\mathcal
{H}_{s}]\\
&=&\mathbb{\hat{E}}[\sum\limits_{i=0}^{k-2}\xi_{t_{i}^{n}}(\langle
M\rangle _{t_{i+1}^{n}}-\langle M\rangle
_{t_{i}^{n}})-2\sum\limits_{i=0}^{k-2}G(\xi_{t_{i}^{n}})( A
_{t_{i+1}^{n}}-A _{t_{i}^{n}})
\\&&+\xi_{t_{k-1}^{n}}(\langle M\rangle
_{s}-\langle M\rangle _{t_{k-1}^{n}})-2G(\xi_{t_{k-1}^{n}})( A
_{s}-A_{t_{k-1}^{n}})\\
&&+\xi_{t_{k-1}^{n}}(\langle M\rangle_{t_{k}^{n}}-\langle M\rangle
_{s})-2G(\xi_{t_{k-1}^{n}})( A _{t_{k}^{n}}-A _{s})|\mathcal
{H}_{s}]\\
&=&X_{s}^{n}+\mathbb{\hat{E}}[\xi_{t_{k-1}^{n}}(\langle
M\rangle_{t_{k}^{n}}-\langle M\rangle
_{s})-2G(\xi_{t_{k-1}^{n}})( A _{t_{k}^{n}}-A _{s})|\mathcal {H}_{s}]\\
&=&X_{s}^{n}+\xi_{t_{k-1}^{n}}^{+}\mathbb{\hat{E}}[\langle M\rangle
_{t_{k}^{n}}-\langle M\rangle _{s}-(A _{t_{k}^{n}}-A _{s})|\mathcal
{H}_{s}]\\
&&+\xi_{t_{k-1}^{n}}^{-}\mathbb{\hat{E}}[-\langle M\rangle
_{t_{k}^{n}}+\langle M\rangle _{s}+\sigma_{0}^{2}(A _{t_{k}^{n}}-A
_{s})|\mathcal{H}_{s}]\\
&=&X_{s}^{n}.
\end{eqnarray*}
For $f\in M_{G,A}^{1}(0,T)$, there exists an $f^{n}\in
M_{G}^{1,0}(0,T)$ such that $$\mathbb{\hat{E}}[\int_0^T
 |f_{s}-f_{s}^{n}|d\langle M\rangle_{s}]\leq\mathbb{\hat{E}}[\int_0^T
 |f_{s}-f_{s}^{n}|dA_{s}]\rightarrow 0,\ \text{as}\ n\rightarrow\infty.$$
 Therefore,
\begin{eqnarray*}
 &&\mathbb{\hat{E}}[|\mathbb{\hat{E}}[X_{t}|\mathcal{H}_{s}]-X_{s}|]\\
&\leq &\mathbb{\hat{E}}[|\mathbb{\hat{E}}[\int_0^t f_{r}d\langle
M\rangle_{r}-2\int_0^t G(f_{r})d A_{r}|\mathcal
{H}_{s}]-\mathbb{\hat{E}}[\int_0^t f_{r}^{n}d\langle
M\rangle_{r}-2\int_0^t G(f_{r}^{n})d A_{r}|\mathcal {F}_{s}]|]
\\&&+\mathbb{\hat{E}}[|\mathbb{\hat{E}}[\int_0^t f_{r}^{n}d\langle
M\rangle_{r}-2\int_0^t G(f_{r}^{n})d A_{r}|\mathcal
{H}_{s}]-\int_0^s f_{r}^{n}d\langle M\rangle_{r}+2\int_0^s
G(f_{r}^{n})d A_{r}|\mathcal {F}_{s}|]
\\&&+\mathbb{\hat{E}}[|\int_0^s f_{r}^{n}d\langle M\rangle_{r}-2\int_0^s
G(f_{r}^{n})d A_{r}-\int_0^s f_{r}d\langle
M\rangle_{r}+2\int_0^s G(f_{r})d A_{r}|]\\
 &\leq &\mathbb{\hat{E}}[|\int_0^t f_{r}d\langle M\rangle_{r}-\int_0^t
f_{r}^{n}d\langle M\rangle_{r}|]+2\mathbb{\hat{E}}[|\int_0^t
G(f_{r})d A_{r}-\int_0^t G(f_{r}^{n})d A_{r}|]
\\&&+\mathbb{\hat{E}}[|\int_0^s f_{r}d\langle M\rangle_{r}-\int_0^s
f_{r}^{n}d\langle M\rangle_{r}|]+2\mathbb{\hat{E}}[|\int_0^s
G(f_{r})d A_{r}-\int_0^s G(f_{r}^{n})d A_{r}|]\\
&\leq & 2\mathbb{\hat{E}}[\int_0^T |f_{s}-f_{s}^{n}|d\langle
M\rangle_{s}]+2\mathbb{\hat{E}}[\int_0^T |f_{s}-f_{s}^{n}|dA_{s}]\\
&\rightarrow & 0,\ \text{as}\ n\rightarrow\infty.
\end{eqnarray*}
Thus,
\begin{eqnarray*}
\mathbb{\hat{E}}[X_{t}|\mathcal {H}_{s}]=X_{s}, \ q.s. \ \text{for
all}\ 0\leq s \leq t\leq T.
\end{eqnarray*}
The proof is complete.
\end{proof}

We can easily get the following corollary, which was established by
Peng \cite{Peng:2007}.
\begin{corollary}
If  $f\in M_{G}^{1}(0,T)$, then  the  process
\begin{eqnarray*}
\Big\{\int_0^t f_{s}d\langle B\rangle_{s}-2\int_0^t G(f_{s})ds , \
t\in[0,T]\Big\} \ \text{is a G-martingale.}
\end{eqnarray*}

\end{corollary}

\begin{remark}
With respect to a linear expectation, if $X$ is a continuous
martingale  with finite  variation, then $X$ is a constant. But it
is not true in G-stochastic analysis. We give an example as follows.
$\{\langle B\rangle_{t}-t\}_{t\geq 0}$ is a continuous G-martingale
with finite
 variation.  But $\{\langle B\rangle_{t}-t\}_{t\geq 0}$ is not a constant. It is a decreasing
  stochastic process.

\end{remark}

\section{Representation of G-martingales as stochastic integrals with respect to  G-Brownian motion}
In this section, we investigate a representation of G-martingales as
stochastic integrals with respect to  G-Brownian motion. The result
of this section will play an important role in the study of
stochastic differential equations driven by G-Brownian motion.

The following martingale characterization of G-Brownian motion was
established by Xu \cite{X:2009}.
\begin{lemma}\label{le}
A process $M\in M_{G}^{2}(0,T)$ is a G-Brownian motion with a
parameter $0<\sigma_{0}\leq
 1$ if
\begin{enumerate} {}
\item[(i)] $M$ is a symmetric G-martingale;
\item[(ii)] For any $t\geq 0$, $M^{2}_{t}-t$ is a  G-martingale;
\item[(iii)] For any $t\geq 0$, $\mathbb{\hat{E}}[-M^{2}_{t}]=-\sigma_{0}^{2}t$;
\item[(iv)] $M$ is continuous, which means for every $\omega\in
\Omega$, $M(t,\omega)$ is continuous.
\end{enumerate}
\end{lemma}

\begin{remark}\label{re2}
It can be easily show that we do not need the assumption $M\in
M_{G}^{2}(0,T)$ in Lemma \ref{le} in our framework. Indeed, one can
use the argument (\ref{e1}) in which $A_{t}$ is replaced by $t$.
\end{remark}

 The following representation of G-martingales as stochastic integrals with respect to  G-Brownian
 motion is the main result in this section.

\begin{theorem}
 Let $0<\sigma_{0}\leq 1$ and  $f\in M_{G}^{2}(0,T)$ be such that $\mathbb{\hat{E}}[\int_0^T|f_{s}|^{4}ds]<\infty$. Moreover,
if there exists a constant $C$ (small enough) such that  $0<C\leq
|f|$, then the following statements $(i)$ is equivalent to $(ii)$.
\begin{enumerate}
\item[(i)] $M$ is a symmetric G-martingale and
  $\Big\{M^{2}_{t}-\int_{0}^{t}f_{s}^{2}ds\Big\}_{t\in[0,T]}$  and
   $\Big\{-M^{2}_{t}+ \sigma_{0}^{2}\int_{0}^{t}f_{s}^{2}ds\Big\}_{t\in[0,T]} $
   are  G-martingales;
\item[(ii)] There exists a G-Brownian motion $B$ such that  $M_{t}=\int_0^t f_{s}dB_{s}$, for all $t\in [0,T]$.
\end{enumerate}
Recall that $G(\alpha)=\frac{1}{2}(\alpha^{+}-\sigma_{0}^{2}\alpha^{-}), \quad \alpha \in \mathbb{R}$.
\end{theorem}

\begin{proof}
We first prove $(i)$  $\Rightarrow(ii)$.

For all $0\leq t \leq T$,   we use
$\pi^{n}=\{0=t_{0}^{n}<t_{1}^{n}\cdots<t_{n}^{n}=t\}$ to denote a
partition of $[0,t]$ such that
  $\max\{t_{i+1}^{n}-t_{i}^{n}, 0\leq i\leq n-1\}\rightarrow 0,$ as $n\rightarrow\infty.$

Since  $f\in M_{G}^{2}(0,T)$, there exists a $f^{n}$ of the form $
f^{n}_{s}=\sum\limits_{j=0}^{n-1}\xi_{t_{j}^{n}}I_{[t_{j}^{n},t_{j+1}^{n})}(s),$
$\xi_{t_{j}^{n}}\in L_{G}^{2}(\mathcal {F}_{t_{j}^{n}}),$ such that
\begin{eqnarray}\label{eq1}
\mathbb{\hat{E}}[\int_0^t| f_{s}-f_{s}^{n}|^{2}ds]\rightarrow 0,\
\text{as}\ n\rightarrow\infty.
\end{eqnarray}
Thanks to  $0<C\leq |f|$, we have
\begin{eqnarray*}\label{}
&&\mathbb{\hat{E}}[|\int_0^t \Big|\frac{1 }{f_{s}}-
\frac{1}{(f_{s}^{n})}\Big|^{2}f_{s}^{2}ds]\\
&\leq & \frac{1}{C^{2}}\mathbb{\hat{E}}[\int_0^t|
f_{s}^{2}-(f_{s}^{n})^{2}|ds]\rightarrow 0,\ \text{as}\
n\rightarrow\infty.
\end{eqnarray*}
Let
 $$X_{t}:=\int_0^t \frac{dM_{s}}{f_{s}}.$$
 Then by Proposition \ref{pr1} and Theorem \ref{th2}, we
know that $X$  is a symmetric G-martingale and continuous. Now we
prove that $\{X^{2}_{t}-t\}_{t\in[0,T]}$ is a G-martingale and
$\mathbb{\hat{E}}[-X^{2}_{t}]=-\sigma_{0}^{2}t$, for all $t\in
[0,T]$.

From the assumptions of $f$ as well as Proposition \ref{pr4} it
follows that the quadratic variation of $M$  exists.

By inequality (\ref{eq1}), we have, for any $\varepsilon>0$,
\begin{eqnarray*}
\mathbb{\hat{E}}[\int_0^t| f_{s}^{2}-(f_{s}^{n})^{2}|ds]&\leq&
(1+\frac{1}{\varepsilon})\mathbb{\hat{E}}[\int_0^t|
f_{s}-f_{s}^{n}|^{2}ds]+ \varepsilon \mathbb{\hat{E}}[\int_0^t|
f_{s}|^{2}ds]\\&\rightarrow& \varepsilon \mathbb{\hat{E}}[\int_0^t|
f_{s}|^{2}ds],\ \text{as}\ n\rightarrow\infty.
\end{eqnarray*}
Therefore,
\begin{eqnarray*}
\mathbb{\hat{E}}[\int_0^t| f_{s}^{2}-(f_{s}^{n})^{2}|ds]\rightarrow
0,\ \text{as}\ n\rightarrow\infty,
\end{eqnarray*}
and
\begin{eqnarray}\label{eq2}
&&\mathbb{\hat{E}}[|\int_0^t \frac{d\langle
M\rangle_{s}}{f_{s}^{2}}-\int_0^t \frac{d\langle
M\rangle_{s}}{(f_{s}^{n})^{2}}|]\nonumber\\
&\leq & \frac{1}{C^{2}}\mathbb{\hat{E}}[\int_0^t|
f_{s}^{2}-(f_{s}^{n})^{2}|ds]\rightarrow 0,\ \text{as}\
n\rightarrow\infty.
\end{eqnarray}
From the subadditivity  of the G-expectation it follows that
\begin{eqnarray*}
\mathbb{\hat{E}}[-X_{t}^{2}]&=&\mathbb{\hat{E}}[-(\int_0^t \frac{d
M_{r}}{f_{r}})^{2}]\\
&\leq &\mathbb{\hat{E}}[-(\int_0^t \frac{d M_{r}}{f_{r}^{n}})^{2}]+
\mathbb{\hat{E}}[-(\int_0^t \frac{d
M_{r}}{f_{r}}-\int_0^t \frac{d M_{r}}{f_{r}^{n}})^{2}]\\
&&+2\mathbb{\hat{E}}[-(\int_0^t \frac{d M_{r}}{f_{r}^{n}})(\int_0^t
\frac{d M_{r}}{f_{r}}-\int_0^t \frac{d M_{r}}{f_{r}^{n}})]\\
\end{eqnarray*}
By the inequality $2ab\leq \varepsilon a^{2}+
\dfrac{1}{\varepsilon}b^{2}$, for all $a,b\in \mathbb{R}$ and
$\varepsilon>0,$ we have

\begin{eqnarray*}
\mathbb{\hat{E}}[-X_{t}^{2}]
 &\leq &\mathbb{\hat{E}}[-(\int_0^t \frac{d M_{r}}{f_{r}^{n}})^{2}]+
\frac{1}{\varepsilon}\mathbb{\hat{E}}[(\int_0^t \frac{d
M_{r}}{f_{r}}-\int_0^t \frac{d M_{r}}{f_{r}^{n}})^{2}]+\varepsilon
\mathbb{\hat{E}}[(\int_0^t \frac{d
M_{r}}{f_{r}^{n}})^{2}]\\
&=&\mathbb{\hat{E}}[-\sum\limits_{i=0}^{n-1}\xi_{t_{i}^{n}}^{-2}(
M^{2} _{t_{i+1}^{n}}-M^{2}
_{t_{i}^{n}})]+\frac{1}{\varepsilon}\mathbb{\hat{E}}[(\int_0^t
\frac{d M_{r}}{f_{r}}-\int_0^t \frac{d
M_{r}}{f_{r}^{n}})^{2}]+\varepsilon \mathbb{\hat{E}}[\int_0^t
\frac{d\langle M\rangle_{s}}{(f_{s}^{n})^{2}}]\\
&\leq&\mathbb{\hat{E}}[-\sum\limits_{i=0}^{n-1}\xi_{t_{i}^{n}}^{-2}
(M^{2}_{t_{i+1}^{n}}-M^{2}_{t_{i}^{n}}-\sigma_{0}^{2}\int_{t_{i}^{n}}^{t_{i+1}^{n}}f_{s}^{2}ds)]
+\sigma_{0}^{2}\mathbb{\hat{E}}[-\sum\limits_{i=0}^{n-1}\xi_{t_{i}^{n}}^{-2}\int_{t_{i}^{n}}^{t_{i+1}^{n}}f_{s}^{2}ds]\\
&&+\frac{1}{\varepsilon}\mathbb{\hat{E}}[(\int_0^t \frac{d
M_{r}}{f_{r}}-\int_0^t \frac{d M_{r}}{f_{r}^{n}})^{2}]+\varepsilon
\mathbb{\hat{E}}[\int_0^t\frac{d\langle M\rangle_{s}}{(f_{s}^{n})^{2}}]\\
&\leq&\mathbb{\hat{E}}[-\sum\limits_{i=0}^{n-1}\xi_{t_{i}^{n}}^{-2}
(M^{2}_{t_{i+1}^{n}}-M^{2}_{t_{i}^{n}}-\sigma_{0}^{2}\int_{t_{i}^{n}}^{t_{i+1}^{n}}f_{s}^{2}ds)]
+\frac{1}{\varepsilon}\mathbb{\hat{E}}[(\int_0^t \frac{d
M_{r}}{f_{r}}-\int_0^t \frac{dM_{r}}{f_{r}^{n}})^{2}]\\
&&+\varepsilon\mathbb{\hat{E}}[\int_0^t\frac{d\langle
M\rangle_{s}}{(f_{s}^{n})^{2}}]
+\sigma_{0}^{2}\mathbb{\hat{E}}[-\sum\limits_{i=0}^{n-1}\xi_{t_{i}^{n}}^{-2}\int_{t_{i}^{n}}^{t_{i+1}^{n}}f_{s}^{2}ds
\\&&+\sum\limits_{i=0}^{n-1}\xi_{t_{i}^{n}}^{-2}\int_{t_{i}^{n}}^{t_{i+1}^{n}}(f_{s}^{n})^{2}ds]+
\sigma_{0}^{2}\mathbb{\hat{E}}[-\sum\limits_{i=0}^{n-1}\xi_{t_{i}^{n}}^{-2}
\int_{t_{i}^{n}}^{t_{i+1}^{n}}(f_{s}^{n})^{2}ds].
\end{eqnarray*}
Thanks to  $(i)$, we obtain
\begin{eqnarray*}
\mathbb{\hat{E}}[-\sum\limits_{i=0}^{n-1}\xi_{t_{i}^{n}}^{-2}
(M^{2}_{t_{i+1}^{n}}-M^{2}_{t_{i}^{n}}-\sigma_{0}^{2}\int_{t_{i}^{n}}^{t_{i+1}^{n}}f_{s}^{2}ds)]=0.
\end{eqnarray*}
Therefore, from inequalities (\ref{eq1}) and (\ref{eq2}) we have
\begin{eqnarray*}
\mathbb{\hat{E}}[-X_{t}^{2}]
&\leq&\sigma_{0}^{2}\mathbb{\hat{E}}[-\sum\limits_{i=0}^{n-1}\xi_{t_{i}^{n}}^{-2}(\int_{t_{i}^{n}}^{t_{i+1}^{n}}f_{s}^{2}ds
-\int_{t_{i}^{n}}^{t_{i+1}^{n}}(f_{s}^{n})^{2}ds)]-\sigma_{0}^{2}t\\
&&+\frac{1}{\varepsilon}\mathbb{\hat{E}}[(\int_0^t \frac{d
M_{r}}{f_{r}}-\int_0^t \frac{d M_{r}}{f_{r}^{n}})^{2}]+\varepsilon
\mathbb{\hat{E}}[\int_0^t\frac{d\langle
M\rangle_{s}}{(f_{s}^{n})^{2}}]\\
 &\leq&\frac{1}{\varepsilon}\mathbb{\hat{E}}[(\int_0^t \frac{d
M_{r}}{f_{r}}-\int_0^t \frac{d
M_{r}}{f_{r}^{n}})^{2}]+\frac{\varepsilon}{C^{2}}
\mathbb{\hat{E}}[\langle M\rangle_{T}]
+\frac{\sigma_{0}^{2}}{C^{2}}\mathbb{\hat{E}}[\int_0^t| f_{s}^{2}-(f_{s}^{n})^{2}|ds]-\sigma_{0}^{2}t\\
&\rightarrow& -\sigma^{2}_{0}t+\frac{\varepsilon}{C^{2}}
\mathbb{\hat{E}}[\langle M\rangle_{T}],\ \text{as}\
n\rightarrow\infty.
\end{eqnarray*}
On the other hand, from the subadditivity  of the G-expectation it
follows that

\begin{eqnarray*}
\mathbb{\hat{E}}[-X_{t}^{2}] &\geq &\mathbb{\hat{E}}[-(\int_0^t
\frac{d M_{r}}{f_{r}^{n}})^{2}]- \mathbb{\hat{E}}[(\int_0^t \frac{d
M_{r}}{f_{r}}-\int_0^t \frac{d M_{r}}{f_{r}^{n}})^{2}]\\
&&-2\mathbb{\hat{E}}[(\int_0^t \frac{d M_{r}}{f_{r}^{n}})(\int_0^t
\frac{d M_{r}}{f_{r}}-\int_0^t \frac{d M_{r}}{f_{r}^{n}})]\\
 &\geq &\mathbb{\hat{E}}[-(\int_0^t \frac{d M_{r}}{f_{r}^{n}})^{2}]-
(1+\frac{1}{\varepsilon})\mathbb{\hat{E}}[(\int_0^t \frac{d
M_{r}}{f_{r}}-\int_0^t \frac{d M_{r}}{f_{r}^{n}})^{2}]-\varepsilon
\mathbb{\hat{E}}[(\int_0^t \frac{d
M_{r}}{f_{r}^{n}})^{2}]\\
&=&\mathbb{\hat{E}}[-\sum\limits_{i=0}^{n-1}\xi_{t_{i}^{n}}^{-2}(
M^{2} _{t_{i+1}^{n}}-M^{2} _{t_{i}^{n}})]-
(1+\frac{1}{\varepsilon})\mathbb{\hat{E}}[(\int_0^t \frac{d
M_{r}}{f_{r}}-\int_0^t \frac{d M_{r}}{f_{r}^{n}})^{2}]-\varepsilon
\mathbb{\hat{E}}[(\int_0^t \frac{d
M_{r}}{f_{r}^{n}})^{2}]\\
&\geq&\mathbb{\hat{E}}[-\sum\limits_{i=0}^{n-1}\xi_{t_{i}^{n}}^{-2}
(M^{2}_{t_{i+1}^{n}}-M^{2}_{t_{i}^{n}}-\sigma_{0}^{2}\int_{t_{i}^{n}}^{t_{i+1}^{n}}f_{s}^{2}ds)]
-\sigma_{0}^{2}\mathbb{\hat{E}}[\sum\limits_{i=0}^{n-1}\xi_{t_{i}^{n}}^{-2}\int_{t_{i}^{n}}^{t_{i+1}^{n}}f_{s}^{2}ds]\\
&&- (1+\frac{1}{\varepsilon})\mathbb{\hat{E}}[(\int_0^t \frac{d
M_{r}}{f_{r}}-\int_0^t \frac{d M_{r}}{f_{r}^{n}})^{2}]-\varepsilon
\mathbb{\hat{E}}[(\int_0^t \frac{d M_{r}}{f_{r}^{n}})^{2}].
\end{eqnarray*}
From  $(i)$ and the inequalities (\ref{eq1}) and (\ref{eq2}) again
it follows that
\begin{eqnarray*} \mathbb{\hat{E}}[-X_{t}^{2}]
&\geq& - (1+\frac{1}{\varepsilon})\mathbb{\hat{E}}[(\int_0^t \frac{d
M_{r}}{f_{r}}-\int_0^t \frac{d M_{r}}{f_{r}^{n}})^{2}]-\varepsilon
\mathbb{\hat{E}}[(\int_0^t \frac{dM_{r}}{f_{r}^{n}})^{2}]\\
&&-\sigma_{0}^{2}\mathbb{\hat{E}}[\sum\limits_{i=0}^{n-1}\xi_{t_{i}^{n}}^{-2}\int_{t_{i}^{n}}^{t_{i+1}^{n}}f_{s}^{2}ds
-\sum\limits_{i=0}^{n-1}\xi_{t_{i}^{n}}^{-2}\int_{t_{i}^{n}}^{t_{i+1}^{n}}(f_{s}^{n})^{2}ds]
-\sigma_{0}^{2}\mathbb{\hat{E}}[\sum\limits_{i=0}^{n-1}\xi_{t_{i}^{n}}^{-2}
\int_{t_{i}^{n}}^{t_{i+1}^{n}}(f_{s}^{n})^{2}ds]\\
&=&- (1+\frac{1}{\varepsilon})\mathbb{\hat{E}}[(\int_0^t \frac{d
M_{r}}{f_{r}}-\int_0^t \frac{d M_{r}}{f_{r}^{n}})^{2}]-\varepsilon
\mathbb{\hat{E}}[(\int_0^t \frac{dM_{r}}{f_{r}^{n}})^{2}]\\
&&-\sigma_{0}^{2}\mathbb{\hat{E}}[\sum\limits_{i=0}^{n-1}\xi_{t_{i}^{n}}^{-2}
\Big(\int_{t_{i}^{n}}^{t_{i+1}^{n}}f_{s}^{2}ds
-\int_{t_{i}^{n}}^{t_{i+1}^{n}}(f_{s}^{n})^{2}ds\Big)]-\sigma_{0}^{2}t\\
 &\geq&- (1+\frac{1}{\varepsilon})\mathbb{\hat{E}}[(\int_0^t \frac{d
M_{r}}{f_{r}}-\int_0^t \frac{d
M_{r}}{f_{r}^{n}})^{2}]-\frac{\varepsilon}{C^{2}}
\mathbb{\hat{E}}[\langle M\rangle_{T}]
-\frac{\sigma_{0}^{2}}{C^{2}}\mathbb{\hat{E}}[\int_0^t| f_{s}^{2}-(f_{s}^{n})^{2}|ds]-\sigma_{0}^{2}t\\
&\rightarrow& -\sigma^{2}_{0}t-\frac{\varepsilon}{C^{2}}
\mathbb{\hat{E}}[\langle M\rangle_{T}],\ \text{as}\
n\rightarrow\infty.
\end{eqnarray*}
Thus, $\mathbb{\hat{E}}[-X_{t}^{2}]=-\sigma^{2}_{0}t,$ for all
$t\in[0,T].$

Now we prove that $\{X^{2}_{t}-t\}_{t\in[0,T]}$ is a G-martingale.
Let $0\leq s \leq t\leq T$. Then
\begin{eqnarray}\label{eq7}
\mathbb{\hat{E}}[X_{t}^{2}-X_{s}^{2}|\mathcal
{H}_{s}]&=&\mathbb{\hat{E}}[(\int_s^t \frac{d
M_{r}}{f_{r}})^{2}+2(\int_s^t \frac{d M_{r}}{f_{r}})(\int_0^s
\frac{d M_{r}}{f_{r}})|\mathcal
{H}_{s}]\nonumber\\
&=&\mathbb{\hat{E}}[(\int_s^t \frac{d M_{r}}{f_{r}})^{2}|\mathcal
{H}_{s}]=\mathbb{\hat{E}}[\int_s^t \frac{d \langle
M\rangle_{r}}{f_{r}^{2}}|\mathcal {H}_{s}]\nonumber\\&\leq&t-s.
\end{eqnarray}
On the other hand, from the subadditivity of G-expectation again it
follows that

\begin{eqnarray*}
\mathbb{\hat{E}}[X_{t}^{2}-X_{s}^{2}|\mathcal {H}_{s}] &\geq
&\mathbb{\hat{E}}[(\int_s^t \frac{d M_{r}}{f_{r}^{n}})^{2}|\mathcal
{H}_{s}]- \mathbb{\hat{E}}[-(\int_s^t \frac{d M_{r}}{f_{r}}-\int_s^t
\frac{dM_{r}}{f_{r}^{n}})^{2}|\mathcal {H}_{s}]\\
 &&-2\mathbb{\hat{E}}[-(\int_s^t \frac{dM_{r}}{f_{r}^{n}})(\int_s^t
\frac{d M_{r}}{f_{r}}-\int_s^t \frac{d M_{r}}{f_{r}^{n}})|\mathcal {H}_{s}]\\
&\geq &\mathbb{\hat{E}}[(\int_s^t \frac{dM_{r}}{f_{r}^{n}})^{2}|\mathcal {H}_{s}]
 -2\mathbb{\hat{E}}[-(\int_s^t \frac{dM_{r}}{f_{r}^{n}})(\int_s^t
\frac{d M_{r}}{f_{r}}-\int_s^t \frac{d M_{r}}{f_{r}^{n}})|\mathcal {H}_{s}]\\
 &\geq &(1-\varepsilon)\mathbb{\hat{E}}[(\int_s^t \frac{d M_{r}}{f_{r}^{n}})^{2}|\mathcal {H}_{s}]
 -\frac{1}{\varepsilon}\mathbb{\hat{E}}[(\int_s^t \frac{d
M_{r}}{f_{r}}-\int_s^t \frac{d M_{r}}{f_{r}^{n}})^{2}|\mathcal {H}_{s}]\\
&=&(1-\varepsilon)\mathbb{\hat{E}}[\sum\limits_{i=0}^{n-1}\xi_{t_{i}^{n}}^{-2}(
M^{2}_{t_{i+1}^{n}}-M^{2}_{t_{i}^{n}\vee s})|\mathcal
{H}_{s}]-\frac{1}{\varepsilon}\mathbb{\hat{E}}[(\int_s^t \frac{d
M_{r}}{f_{r}}-\int_s^t \frac{d M_{r}}{f_{r}^{n}})^{2}|\mathcal
{H}_{s}].
\end{eqnarray*}
By virtue of  $(i)$, we have
\begin{eqnarray*}
\mathbb{\hat{E}}[X_{t}^{2}-X_{s}^{2}|\mathcal
{H}_{s}]&\geq&(1-\varepsilon)\mathbb{\hat{E}}[\sum\limits_{i=0}^{n-1}\xi_{t_{i}^{n}}^{-2}
(M^{2}_{t_{i+1}^{n}}-M^{2}_{t_{i}^{n}\vee s}-\int_{t_{i}^{n}\vee
s}^{t_{i+1}^{n}\vee s}f_{r}^{2}dr)|\mathcal
{H}_{s}]\\
&&-(1-\varepsilon)\mathbb{\hat{E}}[-\sum\limits_{i=0}^{n-1}\xi_{t_{i}^{n}}^{-2}\int_{t_{i}^{n}\vee
s}^{t_{i+1}^{n}\vee s}f_{r}^{2}dr|\mathcal
{H}_{s}]-\frac{1}{\varepsilon}\mathbb{\hat{E}}[(\int_s^t \frac{d
M_{r}}{f_{r}}-\int_s^t \frac{d M_{r}}{f_{r}^{n}})^{2}|\mathcal
{H}_{s}]\\
&\geq& -(1-\varepsilon)
\mathbb{\hat{E}}[-\sum\limits_{i=0}^{n-1}\xi_{t_{i}^{n}}^{-2}\int_{t_{i}^{n}\vee
s}^{t_{i+1}^{n}\vee s}f_{r}^{2}dr
+\sum\limits_{i=0}^{n-1}\xi_{t_{i}^{n}}^{-2}\int_{t_{i}^{n}\vee
s}^{t_{i+1}^{n}\vee
s}(f_{r}^{n})^{2}dr|\mathcal{H}_{s}]\\
&&-(1-\varepsilon)\mathbb{\hat{E}}[-\sum\limits_{i=0}^{n-1}\xi_{t_{i}^{n}\vee
s}^{-2} \int_{t_{i}^{n}\vee s}^{t_{i+1}^{n}\vee
s}(f_{r}^{n})^{2}dr|\mathcal
{H}_{s}]-\frac{1}{\varepsilon}\mathbb{\hat{E}}[(\int_s^t \frac{d
M_{r}}{f_{r}}-\int_s^t\frac{dM_{r}}{f_{r}^{n}})^{2}|\mathcal
{H}_{s}]\\
 &\geq&-\frac{1-\varepsilon}{C^{2}}\mathbb{\hat{E}}[\int_s^t|
f_{s}^{2}-(f_{r}^{n})^{2}|dr|\mathcal
{H}_{s}]+(1-\varepsilon)(t-s)\\
&&-\frac{1}{\varepsilon}\mathbb{\hat{E}}[(\int_0^t\frac{d
M_{r}}{f_{r}}-\int_s^t \frac{d M_{r}}{f_{r}^{n}})^{2}|\mathcal
{H}_{s}].
\end{eqnarray*}
From inequalities (\ref{eq1}) and (\ref{eq2}) it follows that
\begin{eqnarray*}
\mathbb{\hat{E}}[-\mathbb{\hat{E}}[X_{t}^{2}-X_{s}^{2}|\mathcal
{H}_{s}]]
 &\leq&\frac{1-\varepsilon}{C^{2}}\mathbb{\hat{E}}[\int_s^t|
f_{s}^{2}-(f_{r}^{n})^{2}|dr]-(1-\varepsilon)(t-s)\\
&&+\frac{1}{\varepsilon}\mathbb{\hat{E}}[(\int_0^t\frac{d
M_{r}}{f_{r}}-\int_s^t \frac{d M_{r}}{f_{r}^{n}})^{2}]\\
&\rightarrow& -(1-\varepsilon)(t-s),\ \text{as}\ n\rightarrow\infty.
\end{eqnarray*}
Therefore, \begin{eqnarray*}
\mathbb{\hat{E}}[-\mathbb{\hat{E}}[X_{t}^{2}-X_{s}^{2}|\mathcal
{H}_{s}]+(t-s)]&\leq& 0.
\end{eqnarray*}
The above inequality and (\ref{eq7}) yields
\begin{eqnarray*}
\mathbb{\hat{E}}[X_{t}^{2}-t|\mathcal {H}_{s}]=X_{s}^{2}-s, q.s., \
\text{for all} \ \ 0\leq s \leq t \leq T,
\end{eqnarray*}
which means $\{X_{t}^{2}-t\}_{t\in[0,T]}$ is a G-martingale.
Consequently, from Lemma \ref{le} and Remark \ref{re2} we know that $X$ is a G-Brownian
motion with a parameter $\sigma_{0}$.

Now we  prove $(ii)\Rightarrow(i)$.

 From Peng \cite{Peng:2007} we know  that $M$ is a symmetric
 G-martingale. Now put
\begin{eqnarray*}
Y_{t}:=M^{2}_{t}-\int_{0}^{t}f_{s}^{2}ds, \text{for all} \
t\in[0,T].
\end{eqnarray*}

We use $\pi=\{0=t_{0}^{n}<t_{1}^{n}\cdots<t_{n}^{n}=T\}$ to denote a
partition of $[0,T]$ such that
  $\max\{t_{i+1}^{n}-t_{i}^{n}, 0\leq i\leq n-1\}\rightarrow 0,$ as $n\rightarrow\infty.$

For $f\in M_{G}^{2}(0,T)$, there exists an $f^{n}$ of the form $
f_{t}^{n}=\sum\limits_{j=0}^{n-1}\xi_{t_{j}}I_{[t_{j}^{n},t_{j+1}^{n})}(t),$
 where $\xi_{t_{j}}\in L_{G}^{2}(\mathcal {F}_{t_{j}^{n}}),
 0\leq i\leq n-1$, such that
 \begin{eqnarray}\label{eq3}
\mathbb{\hat{E}}[\int_0^T| f_{s}-f_{s}^{n}|^{2}ds]\rightarrow 0,\
\text{as}\ n\rightarrow\infty.
\end{eqnarray}
Then we have
 \begin{eqnarray}\label{eq4}
\mathbb{\hat{E}}[|\int_0^T f_{s}^{2}ds-\int_0^T
(f_{s}^{n})^{2}ds|]\rightarrow 0,\ \text{as}\ n\rightarrow\infty,
\end{eqnarray}
and
 \begin{eqnarray}\label{eq5}
\mathbb{\hat{E}}[|(\int_0^T f_{s}dB_{s})^{2}-(\int_0^T
f_{s}^{n}dB_{s})^{2}|ds]\rightarrow 0,\ \text{as}\
n\rightarrow\infty.
\end{eqnarray}
 Let
\begin{eqnarray*}
Y_{t}^{n}:=(\int_0^tf_{r}^{n}dB_{r})^{2}-\int_{0}^{t}(f_{r}^{n})^{2}dr,
\text{for all} \ t\in[0,T].
\end{eqnarray*}

Let $0\leq s \leq t \leq T$.  We suppose that $ s = t_{k}^{n} \leq
\cdots \leq t_{k+l-1}^{n} \leq t  \leq t_{k+l}^{n} $, for some $k,
l=1,\cdots, n-1$ such that $k+l\leq n$.  Thus,
\begin{eqnarray*}
\mathbb{\hat{E}}[Y_{t}^{n}|\mathcal {H}_{t_{k+l-1}^{n}}]
&=&\mathbb{\hat{E}}[(\int_0^tf_{r}^{n}dB_{r})^{2}-\int_{0}^{t}(f_{r}^{n})^{2}dr|\mathcal
{H}_{t_{k+l-1}^{n}}]\\
&=&\mathbb{\hat{E}}[(\int_0^{t_{k+l-1}^{n}}f_{r}^{n}dB_{r})^{2}-\int_{0}^{t_{k+l-1}^{n}}(f_{r}^{n})^{2}dr
\\&&+(\xi_{t_{k+l-1}^{n}})^{2}(B _{t}-B _{t_{k+l-1}^{n}})^{2}-(\xi_{t_{k+l-1}^{n}})^{2}(
t-t _{t_{k+l-1}^{n}})|\mathcal
{H}_{t_{k}^{n}}]\\
&=&(\int_0^{t_{k+l-1}^{n}}f_{r}^{n}dB_{r})^{2}-\int_{0}^{t_{k+l-1}^{n}}(f_{r}^{n})^{2}dr
\\&&+\mathbb{\hat{E}}[(\xi_{t_{k+l-1}^{n}})^{2}(B _{t}-B _{t_{k+l-1}^{n}})^{2}-(\xi_{t_{k+l-1}^{n}})^{2}(
t-t _{t_{k+l-1}^{n}})|\mathcal
{H}_{t_{k}^{n}}]\\
&=&(\int_0^{t_{k+l-1}^{n}}f_{r}^{n}dB_{r})^{2}-\int_{0}^{t_{k+l-1}^{n}}(f_{r}^{n})^{2}dr\\
&=&Y_{t_{k+l-1}^{n}}.
\end{eqnarray*}
Therefore,
\begin{eqnarray*}
\mathbb{\hat{E}}[Y_{t}^{n}|\mathcal
{H}_{s}]=\mathbb{\hat{E}}[\mathbb{\hat{E}}[Y_{t}^{n}|\mathcal
{H}_{t_{k+l-1}^{n}}]|\mathcal
{H}_{s}]=\mathbb{\hat{E}}[Y_{t_{k+l-1}^{n}}|\mathcal
{H}_{s}]=\cdots=Y_{s}^{n}.
\end{eqnarray*}
Consequently,
\begin{eqnarray*}
\mathbb{\hat{E}}[|\mathbb{\hat{E}}[Y_{t}|\mathcal{H}_{s}]-Y_{s}|]
&\leq &\mathbb{\hat{E}}[|\mathbb{\hat{E}}[Y_{t}|\mathcal
{H}_{s}]-\mathbb{\hat{E}}[Y_{t}^{n}|\mathcal {H}_{s}]|]
\\&&+\mathbb{\hat{E}}[|\mathbb{\hat{E}}[Y_{t}^{n}|\mathcal
{H}_{s}]-Y_{s}^{n}|]+\mathbb{\hat{E}}[|Y_{s}^{n}-Y_{s}|]\\
&\leq &\mathbb{\hat{E}}[|Y_{t}^{n}-Y_{t}|]
+\mathbb{\hat{E}}[|Y_{s}^{n}-Y_{s}|].
\end{eqnarray*}
From  inequalities (\ref{eq4}) and (\ref{eq5}) it follows that
\begin{eqnarray*}
\mathbb{\hat{E}}[|\mathbb{\hat{E}}[Y_{t}|\mathcal{H}_{s}]-Y_{s}|]
&\leq &\mathbb{\hat{E}}[|Y_{t}^{n}-Y_{t}|]
+\mathbb{\hat{E}}[|Y_{s}^{n}-Y_{s}|]\\
&\rightarrow & 0,\ \text{as}\ n\rightarrow\infty.
\end{eqnarray*}
Thus,
\begin{eqnarray*}
\mathbb{\hat{E}}[Y_{t}|\mathcal {H}_{s}]=Y_{s}, \ q.s. \ \text{for
all}\ 0\leq s \leq t\leq T.
\end{eqnarray*}
In the similar argument we can prove that $-
M^{2}+\sigma_{0}^{2}\int_{0}^{\cdot}f_{s}^{2}ds $ is a G-martingale.
The proof is complete.
\end{proof}

\vspace{4mm}

\noindent{\bf Acknowledgements.}

\vspace{4mm} The  author  is very grateful to  Professor Rainer
Buckdahn and Professor Shige Peng for their  helpful discussions and  suggestions.

\end{document}